\documentclass[11pt]{article}
\usepackage[margin=1in]{geometry}
\usepackage[T1]{fontenc}
\usepackage{lmodern,microtype}
\usepackage{amsmath,amssymb,amsthm,mathtools,bm,mathrsfs}
\usepackage{authblk,enumitem,booktabs,array,tabularx}
\usepackage{xcolor}
\usepackage{graphicx,placeins}
\usepackage{tikz}
\usetikzlibrary{arrows.meta,positioning}
\usepackage{hyperref}
\hypersetup{colorlinks=true,linkcolor=blue!55!black,citecolor=blue!55!black,urlcolor=blue!55!black}
\newtheorem{theorem}{Theorem}[section]
\newtheorem{proposition}[theorem]{Proposition}
\newtheorem{corollary}[theorem]{Corollary}
\theoremstyle{definition}

\theoremstyle{remark}
\newtheorem{remark}[theorem]{Remark}
\newcommand{\R}{\mathbb R}
\newcommand{\E}{\mathbb E}
\newcommand{\Harm}{\mathcal H}
\newcommand{\Range}{\operatorname{im}}
\newcommand{\Ker}{\operatorname{ker}}
\newcommand{\Id}{\mathrm{id}}
\newcommand{\ip}[2]{\langle #1,#2\rangle}
\newcommand{\norm}[1]{\lVert #1\rVert}
\newcommand{\HS}{\mathcal L_2}
\newcommand{\dt}{\mathrm d}
\newcommand{\tr}{\operatorname{tr}}
\newcommand{\cS}{\mathscr S}
\newcommand{\constsheaf}{\underline{\R}}
\DeclareMathOperator{\diag}{diag}
\setlist{itemsep=2pt,topsep=4pt}
\numberwithin{equation}{section}
\title{\textbf{Topology-Preserving Hodge--Laplacian Diffusion on Manifolds with Boundary}}
\author[1]{Chen-Chih Lai \thanks{Corresponding author: Chen-Chih.Lai@mtsu.edu}}
\author[2]{Sean Cottrell}
\author[2,3]{Guo-Wei Wei}
\affil[1]{Department of Mathematical Sciences, Middle Tennessee State University, Murfreesboro, TN, USA}
\affil[2]{Department of Mathematics, University of Georgia, Athens, GA, USA}
\affil[3]{Department of Biochemistry and Molecular Biology, University of Georgia, Athens, GA, USA}
\date{}
\begin{document}
\maketitle
\begin{abstract}
{\bf Overview}
We study deterministic and stochastic Hodge--Laplacian diffusion of differential forms on a fixed smooth compact Riemannian manifold with boundary. Absolute and relative boundary conditions select harmonic spaces representing ordinary and relative cohomology. We investigate the preservation properties of these diffusion processes, their compatible cellular discretizations, and their implications for reverse diffusion models. Numerical experiments on an annulus and a solid torus illustrate the theoretical results.

{\bf Mathematical Significance}
Using the de Rham complex and its constant-sheaf interpretation, we distinguish preservation of harmonic coordinates from preservation of a cohomology class. We establish that deterministic heat flow preserves the harmonic projection and exponentially damps its orthogonal complement. For additive stochastic forcing, we prove that pathwise harmonic preservation is equivalent to the absence of harmonic noise. Exact forcing preserves closedness and the initial cohomology class for closed initial data, whereas harmonic-projected forcing need not preserve closedness. We further show that compatible cellular discretizations retain these properties, with an algebraic extension to fixed weighted cellular sheaves. Separating invariant preservation from time-step stability and covariance accuracy, we derive bounds on leakage arising from approximate harmonic projectors.

{\bf AI Relevance}
We develop a finite-dimensional conditional score formulation that confines reverse diffusion to a prescribed harmonic component or, for closed data, a fixed cohomology class. This provides a mathematically grounded framework for incorporating topological constraints into diffusion-based generative AI models on manifolds. 

\end{abstract}


\section{Introduction}\label{sec:intro}
Topological data analysis increasingly studies not only the topology of a base
space, but also data distributed over that space. Sheaf theory provides a
natural local-to-global language for this setting: stalks encode local data,
restriction maps encode compatibility, and sheaf cohomology measures global
consistency and its obstructions \cite{Ghrist2014Sheaves,HansenGhrist2019}.
Persistent sheaf Laplacians provide the associated multiscale framework
\cite{WeiWei2025PSL}, with applications to image representations
\cite{WangWei2026MPSL} and labeled biological networks
\cite{CottrellEtAl2026AD}.
With stalk inner products specified, the cochain complex of a cellular
sheaf carries a Hodge Laplacian whose kernel represents sheaf cohomology.
With standard stalk inner products, the ordinary graph Laplacian is the
degree-zero constant-sheaf case on a graph; higher-degree cellular Hodge
Laplacians arise from the same construction on a cell complex
\cite{HansenGhrist2019}.

At the molecular scale, Hayes et al.\ use localized persistent sheaf Laplacian
spectra to predict protein flexibility through B-factors
\cite{HayesEtAl2025Flexibility}. Ren et al.\ combine persistent sheaf Laplacian
features with protein-transformer embeddings and physicochemical descriptors
in SheafLapNet to predict mutation-induced changes in protein stability and
solubility \cite{RenEtAl2026Stability}. These applications motivate encoding
heterogeneous local information in a sheaf; they are not numerical tests of
the topology-preserving diffusion studied here.

Differential forms furnish the continuous counterpart relevant here. The
assignments $U\mapsto\Omega^q(U)$ form sheaves over a smooth manifold, and the
exterior derivative assembles them into the de Rham complex. Thus the Hodge
Laplacian on forms may be viewed as the analytic Laplacian associated with a
local-to-global sheaf complex. On a smooth compact manifold with boundary,
the absolute and relative boundary realizations select finite-dimensional
harmonic spaces representing ordinary and relative cohomology, respectively.
Boundary conditions are therefore part of the topological data model, not
merely an implementation detail.

Manifold topological deep learning (MTDL) constructs image-derived vector
fields, decomposes them using boundary-aware Hodge operators, and feeds the
resulting components into a convolutional classifier \cite{LiuEtAl2026}.
Combinatorial network-based MTDL (CNMTDL) instead couples the Hodge-derived
representation to attention-based processing on $0$-cell and $2$-cell branches
for image classification \cite{WachiraEtAl2026CNMTDL}. These representation
pipelines motivate our choice of form-valued data. Compatible geometric
discretizations provide the underlying Hodge operators
\cite{ZhaoEtAl2019,SuTongWei2024}.

For multiscale manifold analysis, Su, Tong, and Wei construct persistent de
Rham--Hodge Laplacians using boundary-aware Cartesian-grid operators in an
Eulerian representation \cite{PersistentHodge2024}. Their construction
compares manifolds across a filtration. The present paper instead studies
diffusion time on a fixed realization and the conditions under which selected
harmonic coordinates or a cohomology class remain invariant. We focus on
differential $1$-forms, with statements in general degree where appropriate.
The geometry, metric, coefficient system, and boundary realization are held
fixed during each trajectory. The deterministic Hodge heat flow damps
non-harmonic modes while leaving the selected harmonic component unchanged.
Random perturbations can support data augmentation and forward noising, but
noise with a nonzero harmonic component changes the harmonic coordinates.
Our question is how to constrain the continuous, discrete, and reverse
dynamics to preserve a selected harmonic projection or, for closed data,
a cohomology class.

\paragraph{What is preserved.}
We distinguish three assertions. First, \emph{harmonic-coordinate preservation}
means $P_B\omega_t=P_B\omega_0$, where $P_B$ is the orthogonal projector onto the
kernel of the chosen boundary Hodge Laplacian. This fixes the component of the
form in the harmonic space and is meaningful even for nonclosed forms.
Second, \emph{cohomology preservation} requires that closed initial data
remain closed and that $\omega_t-\omega_0\in\Range d_B$, where $d_B$ is the differential in
the appropriate absolute or relative complex. Third, preservation of the
topology of a decoded image or segmentation requires a separate statement about
the decoder and its output. Our theorems address the first two assertions. The
base manifold's topology is fixed by assumption; harmonic coefficients are
real-valued features of a form, not Betti numbers or a complete description of
an image's topology.

\paragraph{Constant-sheaf scope and sheaf-valued extension.}
The continuum analysis and numerical experiments use the real constant
coefficient system. Sections~\ref{sec:hodge}--\ref{sec:random} work with the
de Rham complex, while Section~\ref{sec:discrete} formulates compatible
cellular discretizations using full or relative cochain complexes.
Applying this description to support-restricted Cartesian operators requires
identifying their retained degrees of freedom and boundary treatment with
the appropriate complex. We also give an algebraic extension to fixed
weighted cellular sheaves. Continuous diffusion with nonconstant coefficient
sheaves requires an additional analytic realization and is discussed
as a future direction in Section~\ref{sec:discussion}.

\paragraph{Relation to existing theory.}
The analytical foundations are classical boundary Hodge theory
\cite{Schwarz1995,HofmannMitreaMonniaux2011}, Hilbert-complex discretization
\cite{ArnoldFalkWinther2010}, and stochastic evolution equations
\cite{DaPratoZabczyk2014}. Hodge heat equations also have an established finite
element theory \cite{ArnoldChen2017}. Spectral cellular sheaf theory supplies the
discrete local-to-global interpretation and the sheaf Laplacian
\cite{HansenGhrist2019}. Weighted boundary Hodge theory and Cartesian-grid
discretizations \cite{SuTongWei2026} provide a natural extension to spatially
varying data weights. Score-based time reversal \cite{SongEtAl2021} supplies the
finite-dimensional generative formulation. The basic semigroup, Hodge
decomposition, sheaf Hodge theorem, and time-reversal identities are not
claimed as new. Our aim is to distinguish harmonic-coordinate preservation
from cohomology preservation in one boundary-consistent setting, carry
these constraints into discrete stochastic dynamics, and separate them
from time-step stability, covariance accuracy, and projector error.

The paper is organized as follows. Section~\ref{sec:sheaf} introduces the
sheaf-theoretic formulation, and Section~\ref{sec:hodge} specifies the
boundary Hodge operators and Hilbert-complex setting.
Section~\ref{sec:det} establishes deterministic preservation and decay,
while Section~\ref{sec:stoch} characterizes harmonic preservation under
stochastic forcing and describes non-harmonic fluctuations.
Section~\ref{sec:random} distinguishes projected and exact perturbations
and proves cohomology preservation for closed initial data under exact
forcing. Section~\ref{sec:discrete} develops compatible cellular
discretizations, structure-preserving time steps, and bounds on leakage
from approximate harmonic projectors. Section~\ref{sec:reverse} formulates
finite-dimensional conditional reverse diffusion with fixed harmonic
coordinates. Section~\ref{sec:numerical} reports deterministic experiments
on a planar annulus and a solid torus, together with a stochastic annulus
comparison of unprotected, harmonic-projected, and exact noise.
Section~\ref{sec:discussion} discusses limitations, further validation,
and extensions beyond the constant-sheaf setting.

\paragraph{Notation and conventions.}
Throughout this manuscript, $q$ denotes differential-form or cochain degree,
$m$ the manifold dimension, $h$ the mesh scale, and $t$ diffusion time.
The upright labels $\mathrm{tan}$ and $\mathrm{nor}$ denote tangential
(absolute) and normal (relative) boundary realizations.
We write $\cS$ for a cellular sheaf on a finite regular cell complex $K$,
$C^q(K;\cS)$ for its cochains, and $d_{\cS}^q$ for its upward coboundary;
the sheaf subscript is suppressed when clear from context. On a manifold,
$d$ is the exterior derivative and $\delta$ the downward Hodge codifferential;
in a cochain complex, $\delta$ denotes the corresponding downward adjoint.
Superscripts on $d$ and $\delta$ indicate their source degree.
The symbols $\star$, ${}^{*}$, and ${}^{\top}$ denote the Hodge star,
an inner-product adjoint, and real matrix transpose, respectively. All
Laplacians use the nonnegative convention. Spaces and differentials outside
the natural degree range are understood to be zero. Boundary realizations
are specified in Section~\ref{sec:hodge}; matrix and mass-coordinate
conventions are introduced in Section~\ref{sec:discrete}.

\section{Sheaf-theoretic and TDA formulation}\label{sec:sheaf}

\subsection{Sheaves as local-to-global data structures}
Let $X$ be a topological space. A sheaf $\cS$ of real vector spaces assigns a
vector space $\cS(U)$ to each open set $U\subset X$ and a restriction map
$\cS(U\supset V):\cS(U)\to\cS(V)$ whenever $V\subset U$. The restrictions satisfy
identity and composition laws, and compatible local sections on an open cover
glue uniquely to a section on the union. The stalk $\cS_x$ is the space of
germs of local sections at $x$, while $\Gamma(X;\cS)=\cS(X)$ is the space of
global sections. This is the local-to-global viewpoint emphasized in applied
sheaf theory \cite{Ghrist2014Sheaves}.

For data analysis, the base space records the geometry or combinatorics of the
sample, while the sheaf records what data live locally and how neighboring
measurements must agree. On a manifold, the constant sheaf $\constsheaf$ assigns to each open set
$U$ the vector space of locally constant real-valued functions on $U$, with
ordinary restriction of functions. Its stalks are one-dimensional, but its
section space on a disconnected open set need not be one-dimensional.
In the cellular model, the constant sheaf assigns $\R$ to each cell and the
identity to each face-restriction map. Nonconstant sheaves permit
varying stalk dimensions and nontrivial compatibility maps, thereby encoding
local frames, partial correspondences, heterogeneous features, or constraints
that cannot be represented by a scalar graph weight alone.

\subsection{Differential forms as a sheaf complex}
On a smooth manifold $M$, the assignment
$U\mapsto\Omega^q(U)$ is a sheaf, denoted $\boldsymbol\Omega^q$. The exterior
derivative commutes with restrictions and yields the de Rham complex of sheaves
\begin{equation}\label{eq:sheaf-derham}
0\longrightarrow\constsheaf\longrightarrow\boldsymbol\Omega^0
 \xrightarrow{d}\boldsymbol\Omega^1\xrightarrow{d}\cdots
 \xrightarrow{d}\boldsymbol\Omega^m\longrightarrow0.
\end{equation}
The augmented sequence \eqref{eq:sheaf-derham} is exact as a sequence of
sheaves and is a fine resolution of $\constsheaf$. The unaugmented complex
of global sections, beginning with $\Omega^0(M)$ in degree zero, is the
ordinary de Rham complex. Taking global sections need not preserve exactness;
its cohomology is identified with the sheaf cohomology of the constant sheaf,
\begin{equation}\label{eq:derham-sheaf-cohomology}
H^q_{\mathrm{dR}}(M)\simeq H^q(M;\constsheaf),
\end{equation}
and the relative complex similarly represents
$H^q(M,\partial M;\constsheaf)$. Thus the harmonic spaces in this paper are
analytic representatives of constant-sheaf cohomology, selected by the metric,
inner product, and boundary realization.

\subsection{Cellular sheaves, cohomology, and spectral diffusion}
Let $K$ be a finite regular cell complex. 
A cellular sheaf $\cS$ assigns a finite-dimensional vector space $\cS(\sigma)$
to each cell and a linear restriction map
$\cS_{\sigma\leq\tau}:\cS(\sigma)\to\cS(\tau)$ to each face relation.
The defining compatibility is
\[
 \cS_{\rho\leq\tau}
 =\cS_{\sigma\leq\tau}\circ\cS_{\rho\leq\sigma}
 \quad(\rho\leq\sigma\leq\tau),\qquad
 \cS_{\sigma\leq\sigma}=\mathrm{id}_{\cS(\sigma)}.
\]
A global section is a compatible assignment $s_\sigma\in\cS(\sigma)$ with
$\cS_{\sigma\leq\tau}s_\sigma=s_\tau$; its space is $\Gamma(K;\cS)$.
The degree-$q$ cochains are
\begin{equation}\label{eq:sheaf-cochains}
C^q(K;\cS)=\bigoplus_{\substack{\sigma\in K\\\dim\sigma=q}}\cS(\sigma).
\end{equation}
A cochain is an arbitrary assignment to cells of a fixed degree; it need not satisfy the compatibility conditions of a global section.

After orientations are chosen, write $d_{\cS}^q:C^q(K;\cS)\to
C^{q+1}(K;\cS)$ for the sheaf coboundary:
\begin{equation}\label{eq:sheaf-coboundary}
(d_{\cS}^q x)_\tau
 =\sum_{\substack{\sigma\leq\tau\\\dim\sigma=q}}
   [\sigma:\tau]\,\cS_{\sigma\leq\tau}x_\sigma,
\qquad d_{\cS}^{q+1}d_{\cS}^q=0.
\end{equation}
For a simplex $\tau=[v_0,\ldots,v_{q+1}]$ and its face obtained by deleting
$v_i$, the induced signed incidence is $[\sigma:\tau]=(-1)^i$.
Orienting simplices by a global vertex ordering follows
\cite[Section~2.2]{WeiWei2025PSL}; for the cubical annulus we instead use its
oriented cellular incidences. Its cohomology and degree-zero interpretation are
\[
 H^q(K;\cS)=\Ker d_{\cS}^q/\Range d_{\cS}^{q-1},
 \qquad H^0(K;\cS)\simeq\Gamma(K;\cS).
\]
Relative sheaf cohomology is obtained by restricting to cochains that vanish
on a boundary subcomplex \cite{Ghrist2014Sheaves,HansenGhrist2019}.

Equip the stalks with inner products and the cochains with their orthogonal
sum, and write $\delta_{\cS}^q=(d_{\cS}^{q-1})^*$ for the downward adjoint.
Following \cite[Section~2.4]{WeiWei2025PSL}, the degree-$q$ sheaf Hodge
Laplacian is
\begin{equation}\label{eq:sheaf-laplacian}
\begin{aligned}
\Delta_{q,\cS}
 &=d_{\cS}^{q-1}\delta_{\cS}^q+\delta_{\cS}^{q+1}d_{\cS}^q\\
 &=d_{\cS}^{q-1}(d_{\cS}^{q-1})^*+(d_{\cS}^q)^*d_{\cS}^q.
\end{aligned}
\end{equation}
This is $d\delta+\delta d$ in degree $q$. The first term is the
\emph{down-Laplacian} and the second is the \emph{up-Laplacian};
the first is absent at degree zero. Finite-dimensional Hodge theory gives
\begin{equation}\label{eq:sheaf-hodge-theorem}
\begin{aligned}
\Ker\Delta_{q,\cS}&\simeq H^q(K;\cS),\\
C^q(K;\cS)&=\Range d_{\cS}^{q-1}\oplus\Ker\Delta_{q,\cS}
 \oplus\Range\delta_{\cS}^{q+1}.
\end{aligned}
\end{equation}
Consequently, the sheaf diffusion
\begin{equation}\label{eq:cellular-sheaf-diffusion}
\dot x=-\Delta_{q,\cS}x
\end{equation}
converges to the orthogonal projection onto the harmonic sheaf cochains.
For $q=0$, the constant sheaf on a graph with standard stalk inner products gives ordinary graph-Laplacian consensus. Other positive stalk inner products give its weighted counterpart. For higher degrees, the constant sheaf yields the corresponding cellular Hodge-Laplacian flow. 
A nonconstant sheaf replaces scalar
agreement by compatibility of local, possibly heterogeneous, data.

The terms \emph{harmonic spectrum} and \emph{non-harmonic spectrum} follow
the PSL literature: they refer to the zero and positive eigenvalues,
respectively \cite{WeiWei2025PSL,CottrellEtAl2026AD,WangWei2026MPSL}.
Harmonic-coordinate preservation below concerns the coefficients of the evolving
state in the zero eigenspace, not changes in the eigenvalues of a fixed operator.

\subsection{The current model as the constant-sheaf specialization}
The smooth continuum state $\omega\in\Omega^1(M)$ is a global section of
$\boldsymbol\Omega^1$, and the differential complex is the de Rham resolution
of $\constsheaf$; its $L^2$ completion supplies the analytic state space below.
For a boundary-conforming cellular model $(K_h,\partial K_h)$, ordinary scalar
cochains are the cochains of the constant cellular sheaf. The full and relative
complexes compute $H^q(K_h;\constsheaf)$ and
$H^q(K_h,\partial K_h;\constsheaf)$, respectively. Section~\ref{sec:discrete}
uses these complexes as the absolute and relative discrete models, with
coboundary matrices $D_B^q$ and mass-weighted Hodge Laplacians.
The additional identification required for support-restricted Cartesian
operators is discussed in Section~\ref{sec:discrete}.

Table~\ref{tab:sheaf-dictionary} records the dictionary used throughout the
paper. Lowercase $d^q$ denotes an abstract differential and uppercase $D^q$
its matrix in specified bases. This follows the operator/matrix distinction
in \cite[Proposition~3.1]{WeiWei2025PSL}.

\begin{table}[htbp]
\centering
\caption{Continuous, cellular-sheaf, and dynamical terminology. Boundary
labels are suppressed in the middle column; use the relative complex when required.}
\label{tab:sheaf-dictionary}
\begin{tabularx}{\textwidth}{>{\raggedright\arraybackslash}p{0.23\textwidth}
>{\raggedright\arraybackslash}p{0.32\textwidth}X}
\toprule
Continuous model & Cellular sheaf model & Interpretation\\
\midrule
$\Omega^q(M)$ & $C^q(K;\cS)$ & forms and cochains assembled by degree\\
$d_B^q$ & $d_{\cS}^q$; matrix $D_{\cS}^q$ & upward differential/coboundary\\
$\delta_B^{q+1}=(d_B^q)^*$ & $\delta_{\cS}^{q+1}=(d_{\cS}^q)^*$ & downward adjoint\\
$\Delta_{q,B}$ & $\Delta_{q,\cS}$ & nonnegative Hodge/sheaf Laplacian\\
$\Ker\Delta_{q,B}$ & $\Ker\Delta_{q,\cS}\simeq H^q(K;\cS)$ & harmonic cohomology representatives\\
$d_B^{q-1}K_t\,\dt W_t$ & $D_{\cS}^{q-1}\mathcal K_h\,\dt W_t$ & exact forcing for closed data\\
\bottomrule
\end{tabularx}
\end{table}

\begin{figure}[htbp]
\centering
\begin{tikzpicture}[
  node distance=0.55cm,
  box/.style={draw,rounded corners,align=center,text width=10.5cm,minimum height=0.9cm,fill=blue!3},
  arr/.style={-{Stealth[length=2.4mm]},thick}
]
\node[box] (data) {Biomedical image, point cloud, or latent features};
\node[box,below=of data] (base) {Choose a base manifold $M$ or cell complex $K_h$ and a sheaf $\cS$ of local feature spaces with restriction maps};
\node[box,below=of base] (complex) {Form the de Rham/sheaf Hilbert complex: $d_{\cS,B}^{q+1}d_{\cS,B}^q=0$};
\node[box,below=of complex] (lap) {Construct the boundary-aware sheaf Hodge Laplacian $\Delta_{q,\cS,B}$ and harmonic projector $P_{\cS,B}$};
\node[box,below=of lap] (dyn) {Forward diffusion/noising and reverse dynamics constrained to preserve the selected harmonic projection or a fixed sheaf-cohomology class};
\node[box,below=of dyn] (out) {Reconstruct or decode, with output topology evaluated separately};
\draw[arr] (data)--(base);
\draw[arr] (base)--(complex);
\draw[arr] (complex)--(lap);
\draw[arr] (lap)--(dyn);
\draw[arr] (dyn)--(out);
\end{tikzpicture}
\caption{Sheaf-theoretic organization of the diffusion framework. The present analysis uses the constant-sheaf/de Rham specialization. Nonconstant coefficients require a compatible cellular complex or an additional continuous analytic realization; decoding is an application-level step.}
\label{fig:sheaf-scheme}
\end{figure}
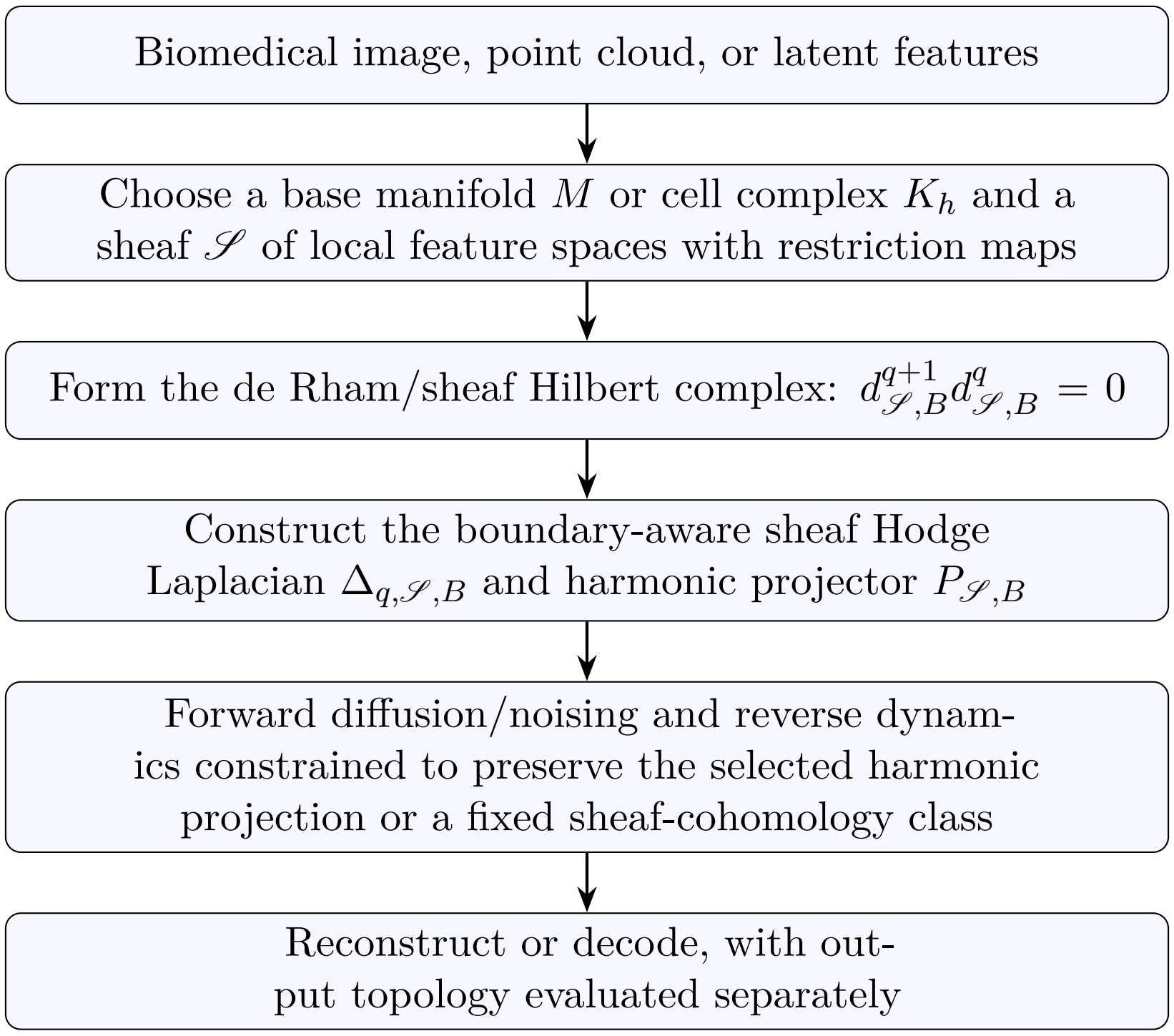

\section{Boundary Hodge theory for the constant sheaf}\label{sec:hodge}
\subsection{Differential forms and boundary conventions}
Let $(M,g)$ be a smooth compact oriented Riemannian manifold of dimension $m$, with
smooth nonempty boundary. All cohomology groups have real coefficients. Following
\cite{ZhaoEtAl2019,PersistentHodge2024,LiuEtAl2026,SuTongWei2026}, write $\Omega^q(M)$ for smooth
$q$-forms, $d:\Omega^q(M)\to\Omega^{q+1}(M)$ for the exterior derivative, and
$\star$ for the Hodge star. The $L^2$ inner product and codifferential are
\begin{equation}\label{eq:continuous-inner-product}
\ip{u}{v}=\int_M u\wedge\star v,\qquad
\delta=(-1)^{m(q-1)+1}\star d\star\quad\text{on }q\text{-forms}.
\end{equation}
The codifferential is the formal adjoint of $d$ in the interior; its Hilbert-space
adjoint domain depends on boundary conditions. 
We write $L^2\Omega^q(M)$ for the space of square-integrable $q$-forms, $H^q_{\mathrm{dR}}(M)$ for de Rham cohomology, and $H^s\Omega^q(M)$ for Sobolev spaces of $q$-forms.
The Hodge Laplacian uses the nonnegative
convention $\Delta_q=d\delta+\delta d$.

A form is \emph{closed} if $d\omega=0$, \emph{exact} if $\omega=d\alpha$,
\emph{coclosed} if $\delta\omega=0$, and \emph{coexact} if $\omega=\delta\gamma$.
In three-dimensional Euclidean space, for a $1$-form associated with a vector field
$v$, $d\omega$ represents $\nabla\times v$ and $\delta\omega$ represents
$-\nabla\cdot v$. Exact $1$-forms are gradient fields and coexact $1$-forms are
curl fields. They are respectively curl-free and divergence-free, but the converse
implications need not hold globally. We therefore use \emph{exact} and
\emph{coexact} for the summands, and \emph{closed} and \emph{coclosed} for the
corresponding differential constraints. The form terminology applies in any dimension.

Let $\iota:\partial M\hookrightarrow M$ be inclusion. The first-order trace spaces
in the cited papers are
\begin{equation}\label{eq:normal-tangential-spaces}
\Omega_{\mathrm{nor}}^q=\{\omega\in\Omega^q(M):\iota^*\omega=0\},\qquad
\Omega_{\mathrm{tan}}^q=\{\omega\in\Omega^q(M):\iota^*(\star\omega)=0\}.
\end{equation}
We use $B\in\{\mathrm{tan},\mathrm{nor}\}$ for the
\emph{tangential} and \emph{normal} realizations, respectively.
For $1$-forms, these mean respectively $v\cdot\mathbf n=0$ and vanishing tangential component of $v$, where $\mathbf n$ is the outward unit normal. For other degrees,
the definitions are the trace conditions in \eqref{eq:normal-tangential-spaces},
not the same vector-field interpretation. The full tangential Laplacian is also
called the \emph{absolute} realization; the full normal Laplacian is the
\emph{relative} realization.

The terms Neumann and Dirichlet in the references extend scalar terminology.
For functions, the full tangential realization gives the scalar Neumann condition,
and the full normal realization gives the scalar Dirichlet condition. On positive-degree
forms they impose the two trace conditions below, not componentwise scalar
Neumann or Dirichlet conditions.

\subsection{Boundary realizations and closed de Rham complexes}
On a smooth boundary, the classical second-order operator domains are
\begin{align}
D(\Delta_{q,\mathrm{tan}})
 &=\{\omega\in H^2\Omega^q(M):
       \iota^*(\star\omega)=0,\ \iota^*(\star d\omega)=0\},
       \label{eq:absolute-bc}\\
D(\Delta_{q,\mathrm{nor}})
 &=\{\omega\in H^2\Omega^q(M):
       \iota^*\omega=0,\ \iota^*(\delta\omega)=0\}.
       \label{eq:relative-bc}
\end{align}
For smooth forms, these domains coincide with the second-order trace
spaces in \cite{SuTongWei2026}. 
For intermediate degrees $1\le q\le m-1$, the trace conditions defining $\Omega_{\mathrm{tan}}^q$ or $\Omega_{\mathrm{nor}}^q$ alone do not specify the full second-order realization. 
At degrees $q=0$ and $q=m$, some of the displayed conditions are automatic.

For smooth $\alpha\in\Omega^{q-1}(M)$ and $\beta\in\Omega^q(M)$, Green's formula is
\begin{equation}\label{eq:green-hodge}
\ip{d\alpha}{\beta}=\ip{\alpha}{\delta\beta}
 +\int_{\partial M}\iota^*(\alpha\wedge\star\beta).
\end{equation}
Section~\ref{sec:det-energy} derives the energy identity from this formula and
checks the cancellation of each boundary term.

To specify weak formulations, define the closed exterior derivatives by
\begin{align}
D(d_{\mathrm{tan}}^q)&=\{\omega\in L^2\Omega^q(M):
           d\omega\in L^2\Omega^{q+1}(M)\text{ distributionally}\},
           \label{eq:absolute-graph-domain}\\
D(d_{\mathrm{nor}}^q)&=\overline{C_c^\infty\Omega^q(M^\circ)}^{\,\norm{\cdot}+\norm{d(\cdot)}}.
           \label{eq:relative-graph-domain}
\end{align}
These give the maximal and minimal de Rham complexes, respectively. The normal
complex imposes zero tangential trace weakly. The tangential realization uses
$d_{\mathrm{tan}}=d_{\max}$ with no essential trace on its domain; its adjoint imposes the zero
normal trace. The label $\mathrm{tan}$ identifies the resulting Laplacian realization;
it does not impose the trace defining $\Omega_{\mathrm{tan}}^q$ on the domain of $d_{\mathrm{tan}}^q$.
Indeed, $d$ preserves $\Omega_{\mathrm{nor}}^q$, whereas $\delta$ preserves $\Omega_{\mathrm{tan}}^q$.
Set $\delta_B^q=(d_B^{q-1})^*$.
The closed nonnegative form
\begin{equation}\label{eq:hodge-form}
\mathfrak a_{q,B}(u,v)=\ip{d_B^qu}{d_B^qv}
             +\ip{\delta_B^qu}{\delta_B^qv},\qquad
V_{q,B}=D(d_B^q)\cap D(\delta_B^q)
\end{equation}
has associated self-adjoint operator $\Delta_{q,B}$ with the domain
\eqref{eq:absolute-bc} or \eqref{eq:relative-bc}
\cite{Schwarz1995,HofmannMitreaMonniaux2011}. These operators have compact
resolvent, and
\begin{equation}\label{eq:exact-hodge-energy}
V_{q,B}=D(\Delta_{q,B}^{1/2}),\qquad
\norm{\Delta_{q,B}^{1/2}u}^2
 =\norm{d_B^qu}^2+\norm{\delta_B^qu}^2.
\end{equation}
This form domain is used in the stochastic energy estimates. It does not require
classical second-order boundary traces of every noise realization.

\subsection{Harmonic fields, harmonic forms, and cohomology}
It is useful to distinguish harmonic fields from Laplacian-harmonic forms,
as in \cite[Section~2]{PersistentHodge2024} and
\cite[Supplementary Information, Section~1.1]{LiuEtAl2026}:
\begin{equation}\label{eq:unrestricted-harmonic-spaces}
\Harm^q(M)=\{\omega\in\Omega^q(M):d\omega=0,\ \delta\omega=0\},\qquad
\Harm_\Delta^q(M)=\{\omega\in\Omega^q(M):\Delta_q\omega=0\}.
\end{equation}
We call elements of $\Harm^q(M)$ \emph{harmonic fields} and elements of
$\Harm_\Delta^q(M)$ \emph{Laplacian-harmonic forms} when no boundary condition
has been specified. On a manifold with boundary the two conditions are not
equivalent: $\Harm^q(M)\subseteq\Harm_\Delta^q(M)$, possibly strictly.
The boundary-conditioned harmonic spaces are
\begin{equation}\label{eq:boundary-harmonic-fields}
\Harm_{\mathrm{tan}}^q=\Harm^q(M)\cap\Omega_{\mathrm{tan}}^q=\Ker\Delta_{q,\mathrm{tan}},\qquad
\Harm_{\mathrm{nor}}^q=\Harm^q(M)\cap\Omega_{\mathrm{nor}}^q=\Ker\Delta_{q,\mathrm{nor}}.
\end{equation}
The last equalities follow from the energy identity; elements of these kernels
are smooth by elliptic regularity. In the rest of the paper, \emph{harmonic
projection} always refers to one of these finite-dimensional kernels.

The closed-range Hilbert-complex decomposition is
\begin{equation}\label{eq:hodge-decomp}
L^2\Omega^q(M)=E_B^q\oplus\Harm_B^q\oplus C_B^q,\qquad
E_B^q=\Range d_B^{q-1},\quad C_B^q=\Range\delta_B^{q+1},
\end{equation}
with mutually orthogonal closed summands. Here exactness and coexactness refer
to the specified operator domains, including their boundary restrictions. Moreover,
\begin{align}\label{eq:closed-decomp}
Z_B^q&:=\Ker d_B^q=E_B^q\oplus\Harm_B^q,\notag\\
\Harm_{\mathrm{tan}}^q&\simeq H_{\mathrm{dR}}^q(M)\simeq H^q(M;\constsheaf),\notag\\
\Harm_{\mathrm{nor}}^q&\simeq H_{\mathrm{dR}}^q(M,\partial M)
\simeq H^q(M,\partial M;\constsheaf).
\end{align}
Write $\beta_j(M)=\dim H_{\mathrm{dR}}^j(M)$ for the Betti numbers, as in the
cited papers. Poincar\'e--Lefschetz duality gives
\begin{equation}\label{eq:nullity}
\dim\Harm_{\mathrm{tan}}^q=\beta_q(M),\qquad \dim\Harm_{\mathrm{nor}}^q=\beta_{m-q}(M).
\end{equation}
For $1$-forms in three dimensions, these dimensions are $\beta_1$ and $\beta_2$;
in two dimensions both are $\beta_1$. Thus planar nullities alone cannot
distinguish the two boundary implementations.

Let $P_B^q$ be the $L^2$-orthogonal projector onto $\Harm_B^q$ and
$Q_B^q=\Id-P_B^q$. For $1$-forms suppress $q$ and write
$H=L^2\Omega^1(M)$, $\Delta_B$, $\Harm_B$, $P_B$, and $Q_B$.
The first positive eigenvalue is $\lambda_{B,+}>0$, and
\begin{equation}\label{eq:poincare}
\norm{\Delta_B^{1/2}u}^2\ge\lambda_{B,+}\norm{u}^2,
\qquad u\in V_{1,B}\cap Q_BH.
\end{equation}
If the harmonic space is trivial, harmonic preservation is vacuous and $Q_B=\Id$.

\subsection{Relation to the three- and five-component decompositions}
The decomposition \eqref{eq:hodge-decomp} for a fixed realization must be
distinguished from the Hodge--Morrey decomposition
\begin{equation}\label{eq:hodge-morrey}
\Omega^q(M)=d\Omega_{\mathrm{nor}}^{q-1}\oplus\delta\Omega_{\mathrm{tan}}^{q+1}\oplus\Harm^q(M)
\end{equation}
used in \cite[Supplementary Information, equation (10)]{LiuEtAl2026} and
\cite[Section 2.4]{ZhaoEtAl2019}. Its last summand has no imposed boundary trace
and can be infinite-dimensional; for example, gradients of scalar harmonic
functions on a ball in $\R^m$, $m\ge2$, supply such fields in degree one.
It is not generally either
$\Harm_{\mathrm{tan}}^q$ or $\Harm_{\mathrm{nor}}^q$.

For smooth compact domains in $\R^3$ and $1$-forms, the five-component
decomposition in \cite{ZhaoEtAl2019} is
\begin{equation}\label{eq:five-component-comparison}
\Omega^1(M)=d\Omega_{\mathrm{nor}}^0\oplus\delta\Omega_{\mathrm{tan}}^2
 \oplus\Harm_{\mathrm{tan}}^1\oplus\Harm_{\mathrm{nor}}^1
 \oplus\bigl(d\Omega^0(M)\cap\delta\Omega^2(M)\bigr).
\end{equation}
All five summands are orthogonal in this Euclidean-domain setting. In vector
terminology the last term is the \emph{curl-gradient component}: it is both a
gradient and a curl, and is therefore closed and coclosed. The orthogonality
of $\Harm_{\mathrm{tan}}^q$ and $\Harm_{\mathrm{nor}}^q$ should not be assumed on a general Riemannian
manifold. The Hodge--Morrey residual in \eqref{eq:hodge-morrey} includes the
last three terms of \eqref{eq:five-component-comparison}; our diffusion fixes
only the projection selected by $B$. For example, a nonzero constant $1$-form
on a ball belongs to $\Harm^1(M)$, while both degree-one boundary-conditioned
kernels are zero. Its full Hodge--Morrey residual is consequently not preserved
by either boundary-conditioned heat flow.

This distinction also matters computationally. 
The harmonic residual used in an image decomposition need not equal the projection onto the kernel of the boundary-conditioned Laplacian used for diffusion.
Identifying the two requires agreement of the operators, boundary conditions, and inner products. Section~\ref{sec:discrete} states the corresponding algebraic guarantees for a specified compatible complex.

\subsection{Fixed spatial weights}
Following \cite[Sections 3--4]{SuTongWei2026}, let $f:M\to\R$ be a fixed smooth
function and let $\mu_g$ denote Riemannian volume measure. Define
\begin{equation}\label{eq:weighted-hodge}
\ip{u}{v}_f=\int_M\ip{u}{v}_g e^{-f}\,\dt\mu_g,\qquad
\delta_f=e^f\delta e^{-f},\qquad \Delta_f=d\delta_f+\delta_fd.
\end{equation}
The coefficient sheaf and exterior derivative remain unchanged; the weight modifies the Hilbert-space inner product, its adjoint, and the associated harmonic representatives.

Since $M$ is compact, the density $e^{-f}$ is bounded above and below
by positive constants. The tangential traces remain those in \eqref{eq:absolute-bc}; the second
normal trace becomes $\iota^*(\delta_f\omega)=0$. Orthogonality, energy, covariance
adjoints, and the harmonic projector must all use $\ip{\cdot}{\cdot}_f$.
The weighted harmonic fields are closed and $f$-coclosed, meaning
$d\omega=0$ and $\delta_f\omega=0$; weighted coexactness means membership
in the appropriate range of $\delta_f$. In Euclidean vector notation,
$\delta_f\omega=0$ expresses $\nabla\cdot(e^{-f}v)=0$, not necessarily
$\nabla\cdot v=0$.

The weighted harmonic spaces represent the same cohomology groups, with
generally different representatives and projectors. The arguments below
apply to each fixed weighted realization when its inner products and
adjoints are used consistently. With the convention
\eqref{eq:weighted-hodge}, $\delta_f=\delta+\iota_{\nabla f}$, where
$\iota_{\nabla f}$ denotes contraction. Other boundary realizations are
discussed in Section~\ref{sec:other-bc}.

\section{Deterministic Hodge--Laplacian diffusion}\label{sec:det}
Fix one of the two boundary realizations from Section~\ref{sec:hodge}, with
$B\in\{\mathrm{tan},\mathrm{nor}\}$. The manifold, metric, and boundary
realization remain fixed.
On $H=L^2\Omega^1(M)$, consider
\begin{equation}\label{eq:hodge-heat}
\partial_t\omega+\nu\Delta_B\omega=0,\qquad
\omega(0)=\omega_0,\qquad \nu>0.
\end{equation}
Since $\Delta_B$ is nonnegative and self-adjoint, the spectral theorem defines
the contraction semigroup
\begin{equation}\label{eq:det-semigroup}
S_B(t)=e^{-\nu t\Delta_B},\qquad \omega(t)=S_B(t)\omega_0.
\end{equation}
The infinitesimal generator of this semigroup is $-\nu\Delta_B$.
Equation~\eqref{eq:det-semigroup} is the unique mild solution for every $\omega_0\in H$, with
$\omega\in C([0,\infty);H)$. For $t>0$, semigroup smoothing gives
$\omega(t)\in D(\Delta_B)$ and the equation holds strongly in $H$.
The initial form need not satisfy classical boundary traces: those traces are
encoded by the generator and hold at positive times in its operator-domain
sense. If $\omega_0\in D(\Delta_B)$, the solution is also strong at the initial
time, with a right derivative there.

\subsection{Harmonic invariance}
Let $\Harm_B=\Ker\Delta_B$, let $P_B$ be its $L^2$-orthogonal projector, and
write $Q_B=\Id-P_B$. The invariance statement uses only nonnegativity and
self-adjointness; the compactness and spectral gap needed for a quantitative
decay rate enter in the next subsection.

\begin{theorem}[Harmonic projection is invariant]\label{thm:det-invariance}
Let $\Delta_B$ be the nonnegative self-adjoint realization in
\eqref{eq:hodge-heat}. For every $\omega_0\in H$, its mild solution satisfies
\begin{equation}\label{eq:det-invariance}
P_B\omega(t)=P_B\omega_0\qquad\text{for every }t\ge0.
\end{equation}
Consequently, for any fixed orthonormal basis $\{h_\ell\}_{\ell=1}^r$ of
$\Harm_B$, each harmonic coefficient
$c_\ell(t)=\ip{\omega(t)}{h_\ell}$ equals $c_\ell(0)$.
\end{theorem}
\begin{proof}
The spectral projector $P_B$ commutes with $S_B(t)$, and $S_B(t)$ is the
identity on $\Ker\Delta_B$. Hence
\[
P_B\omega(t)=P_BS_B(t)\omega_0
 =S_B(t)P_B\omega_0=P_B\omega_0.
\]
All operators in this identity are bounded on $H$, so no operator-domain
assumption on $\omega_0$ is required. Taking inner products with $h_\ell$
gives the coefficient identities.
\end{proof}

The theorem applies to nonclosed forms and preserves their harmonic
coordinates, not a cohomology class. If $\Harm_B=\{0\}$, the statement is
vacuous and the entire solution lies in the dissipative subspace.

\begin{remark}[Constant-sheaf diffusion]
Theorem~\ref{thm:det-invariance} is the continuum counterpart of
\eqref{eq:cellular-sheaf-diffusion}. The degree-zero constant-sheaf case
recovers the graph consensus described in Section~\ref{sec:sheaf}.
\end{remark}

\subsection{Decay of non-harmonic modes}
For the smooth compact setting of Section~\ref{sec:hodge}, $\Delta_B$ has
compact resolvent. Its kernel is finite dimensional, and the positive spectrum
is discrete. Suppose $Q_BH\ne\{0\}$ and write
\begin{equation}\label{eq:det-positive-gap}
\lambda_{B,+}=\min\{\lambda>0:\lambda\in\operatorname{spec}(\Delta_B)\}.
\end{equation}
Choose an orthonormal basis $\{e_j\}_{j\ge1}$ of $Q_BH$ consisting of
eigenforms, with
$\Delta_Be_j=\lambda_je_j$ and $\lambda_j\ge\lambda_{B,+}>0$.

\begin{theorem}[Spectral decay on the orthogonal complement]\label{thm:decay}
For every $\omega_0\in H$ and every $t\ge0$,
\begin{equation}\label{eq:decay-estimate}
\norm{Q_B\omega(t)}
 \le e^{-\nu\lambda_{B,+}t}\norm{Q_B\omega_0}.
\end{equation}
In particular,
\begin{equation}\label{eq:long-time-limit}
\omega(t)\longrightarrow P_B\omega_0
\quad\text{in }H\text{ as }t\to\infty.
\end{equation}
\end{theorem}
\begin{proof}
Write $Q_B\omega_0=\sum_{j\ge1}a_je_j$, with
$\sum_j|a_j|^2=\norm{Q_B\omega_0}^2$. Spectral calculus and Parseval's
identity give
\[
\norm{Q_B\omega(t)}^2
 =\sum_{j\ge1}e^{-2\nu\lambda_jt}|a_j|^2
 \le e^{-2\nu\lambda_{B,+}t}\norm{Q_B\omega_0}^2.
\]
Taking square roots proves \eqref{eq:decay-estimate}. The limit follows
from $\omega(t)-P_B\omega_0=Q_B\omega(t)$.
\end{proof}

The exponent $\nu\lambda_{B,+}$ is the sharp uniform $H$-norm decay rate:
equality holds for an initial positive component entirely in the first positive
eigenspace. A particular initial form can decay faster if its coefficient in that
eigenspace vanishes. Thus a measured decay rate must be compared with the
smallest positive eigenvalue actually excited by the initial field. The squared
norm decays with exponent $2\nu\lambda_{B,+}$.

\subsection{Energy identity and boundary compatibility}\label{sec:det-energy}
For the two standard boundary realizations, the form identity
\eqref{eq:exact-hodge-energy} in degree one is
\begin{equation}\label{eq:det-form-energy}
\mathfrak a_{1,B}(u,u)
 =\norm{d_B^1u}^2+\norm{\delta_B^1u}^2
 =\norm{\Delta_B^{1/2}u}^2,
\qquad u\in V_{1,B}.
\end{equation}
Unlike $\ip{\Delta_Bu}{u}$, which requires $u\in D(\Delta_B)$,
this identity is defined on the form domain $V_{1,B}$.

For a smooth $1$-form $u$, apply Green's formula
\eqref{eq:green-hodge} first to $\alpha=\delta u$, $\beta=u$, and then to
$\alpha=u$, $\beta=du$. This gives
\begin{align*}
\ip{d\delta u}{u}
 &=\norm{\delta u}^2+\int_{\partial M}\iota^*(\delta u\wedge\star u),\\
\ip{\delta du}{u}
 &=\norm{du}^2-\int_{\partial M}\iota^*(u\wedge\star du).
\end{align*}
Adding the two equalities yields the full boundary energy identity
\begin{equation}\label{eq:energy-with-boundary}
\begin{aligned}
\ip{\Delta_1u}{u}
 &=\norm{du}^2+\norm{\delta u}^2\\
 &\quad+\int_{\partial M}\iota^*
       (\delta u\wedge\star u-u\wedge\star du).
\end{aligned}
\end{equation}
Define its quadratic boundary contribution by
\begin{equation}\label{eq:det-quadratic-boundary}
\mathcal B(u)=\int_{\partial M}\iota^*
 \bigl(\delta u\wedge\star u-u\wedge\star du\bigr).
\end{equation}
Then
$\ip{\Delta_1u}{u}=\norm{du}^2+\norm{\delta u}^2+\mathcal B(u)$.
For $B=\mathrm{tan}$, the traces $\iota^*(\star u)=0$ and
$\iota^*(\star du)=0$ make the two boundary terms vanish separately.
For $B=\mathrm{nor}$, the same conclusion follows from
$\iota^*u=0$ and $\iota^*(\delta u)=0$. Hence
\begin{equation}\label{eq:hodge-energy}
\ip{\Delta_Bu}{u}
 =\norm{d_B^1u}^2+\norm{\delta_B^1u}^2,
\qquad u\in D(\Delta_B).
\end{equation}
The trace computation agrees with the closed-form construction; the latter
extends the identity beyond smooth test fields.

\begin{proposition}[Deterministic energy balance]\label{prop:det-energy}
For $\omega_0\in H$, the mild solution belongs to
$L^2(0,T;V_{1,B})$ for every finite $T$ and satisfies
\begin{equation}\label{eq:deterministic-energy-rigorous}
\norm{\omega(t)}^2
 +2\nu\int_0^t\norm{\Delta_B^{1/2}\omega(s)}^2\,\dt s
 =\norm{\omega_0}^2\qquad(t\ge0).
\end{equation}
For every $t>0$, it also satisfies the differential identity
\begin{equation}\label{eq:det-differential-energy}
\frac12\frac{\dt}{\dt t}\norm{\omega(t)}^2
 =-\nu\bigl(\norm{d_B^1\omega(t)}^2
             +\norm{\delta_B^1\omega(t)}^2\bigr).
\end{equation}
If $\omega_0\in D(\Delta_B)$, this identity holds at $t=0$ with the
right derivative as well.
\end{proposition}
\begin{proof}
First take operator-domain initial data. The strong equation and
\eqref{eq:hodge-energy} give
\[
\frac12\frac{\dt}{\dt t}\norm{\omega(t)}^2
 =\ip{\partial_t\omega(t)}{\omega(t)}
 =-\nu\ip{\Delta_B\omega(t)}{\omega(t)}
 =-\nu\mathfrak a_{1,B}(\omega(t),\omega(t)).
\]
This proves \eqref{eq:det-differential-energy}; integrating proves
\eqref{eq:deterministic-energy-rigorous} in this case.

For general $\omega_0\in H$, write
$\omega_0=P_B\omega_0+\sum_{j\ge1}a_je_j$ as in
Theorem~\ref{thm:decay}. The harmonic component contributes no dissipation.
Tonelli's theorem applies to the nonnegative summands and gives
\[
\begin{aligned}
2\nu\int_0^t\norm{\Delta_B^{1/2}\omega(s)}^2\,\dt s
 &=\sum_{j\ge1}2\nu\lambda_j|a_j|^2
       \int_0^te^{-2\nu\lambda_js}\,\dt s\\
 &=\sum_{j\ge1}|a_j|^2(1-e^{-2\nu\lambda_jt})\\
 &=\norm{\omega_0}^2-\norm{\omega(t)}^2.
\end{aligned}
\]
This proves the integrated identity without requiring an initial derivative or
initial trace. Its finiteness, together with the contraction estimate in $H$,
implies $\omega\in L^2(0,T;V_{1,B})$. 
Finally, for every $t>0$, the solution $\omega(t)=S_B(t)\omega_0$ belongs to $D(\Delta_B)$, so the strong calculation applies and proves the differential identity.
\end{proof}

Applying \eqref{eq:poincare} to $Q_B\omega(t)$ recovers the decay
estimate \eqref{eq:decay-estimate} by Gronwall's inequality.
Letting $t\to\infty$ in the integrated energy balance gives
\begin{equation}\label{eq:total-det-dissipation}
2\nu\int_0^\infty
 \bigl(\norm{d_B^1\omega(s)}^2+\norm{\delta_B^1\omega(s)}^2\bigr)\,\dt s
 =\norm{Q_B\omega_0}^2.
\end{equation}

\subsection{Closedness and the initial cohomology class}
Harmonic invariance is enough to preserve a de Rham class only when closedness
is also preserved. For either standard boundary complex and any degree $q$,
the heat semigroup satisfies
\begin{equation}\label{eq:heat-commutes}
d_B^qS_{q,B}(t)u=S_{q+1,B}(t)d_B^qu,
\qquad u\in D(d_B^q),\qquad
S_{q,B}(t)=e^{-\nu t\Delta_{q,B}}.
\end{equation}
The corresponding resolvent commutation is given in
\cite[Proposition 3.1]{HofmannMitreaMonniaux2011}; closedness of $d_B^q$
and resolvent approximation to the exponential give
\eqref{eq:heat-commutes}. In particular the semigroup preserves both the
closed and exact subspaces of the selected complex.

If $\omega_0\in Z_B^q=\Ker d_B^q$, the Hodge decomposition writes it
uniquely as $\omega_0=h_0+e_0$ with $h_0=P_B^q\omega_0$ and
$e_0\in E_B^q=\Range d_B^{q-1}$. Hence
\[
\omega(t)=h_0+S_{q,B}(t)e_0\in h_0+E_B^q,
\qquad \omega(t)-\omega_0\in E_B^q.
\]
The trajectory therefore remains closed and in its initial absolute class for
$B=\mathrm{tan}$, or its initial relative class for $B=\mathrm{nor}$. For nonsmooth forms, these
classes are understood in the $L^2$ Hilbert complex and identified with the
corresponding de Rham cohomology. Section~\ref{sec:random} extends this affine-space argument to exact
random forcing.

\subsection{Other boundary conditions compatible with the energy identity}
\label{sec:other-boundary}\label{sec:other-bc}
The bulk identity \eqref{eq:hodge-energy} also holds for any nonnegative
self-adjoint realization whose operator-domain traces give $\mathcal B(u)=0$.
Boundary cancellation alone establishes neither self-adjointness nor a
cohomological interpretation of the kernel.

\paragraph{Quadratic boundary energy and the Green pairing.}
For real smooth $q$-forms, define the bilinear boundary expression
\begin{equation}\label{eq:bilinear-boundary}
\mathcal C(u,v)=\int_{\partial M}\iota^*
  \bigl(\delta u\wedge\star v-v\wedge\star du\bigr).
\end{equation}
Integration by parts gives
\[
\ip{\Delta_qu}{v}
 =\ip{du}{dv}+\ip{\delta u}{\delta v}+\mathcal C(u,v),
\qquad \mathcal B(u)=\mathcal C(u,u).
\]
The Green pairing relevant to symmetry is instead the skew expression
\begin{equation}\label{eq:green-skew-pairing}
\mathcal G(u,v)
 :=\ip{\Delta_qu}{v}-\ip{u}{\Delta_qv}
 =\mathcal C(u,v)-\mathcal C(v,u).
\end{equation}
Annihilating $\mathcal G$ yields symmetry; self-adjointness further
requires equality with the adjoint domain, and ellipticity requires its
own justification. Symmetry does not imply $\mathcal C(u,u)=0$: the skew
pairing always vanishes on the diagonal, whereas the energy boundary term
need not vanish.

\paragraph{Mixed tangential/normal conditions.}
Let $\Gamma_{\mathrm{tan}}$ and $\Gamma_{\mathrm{nor}}$ be disjoint relatively open boundary pieces
whose closures cover $\partial M$, with a regular common interface. A sufficient
setting for the following discussion is a smooth boundary partition with a
smooth interface; more generally one uses admissible Lipschitz boundary pairs
\cite{Licht2019,PaulySchomburg2022}. For $1$-forms, our subscripts name
the allowed vector-field direction: $\Gamma_{\mathrm{tan}}$ carries the tangential/absolute realization, and
$\Gamma_{\mathrm{nor}}$ the normal/relative realization. In particular, the vanishing
\emph{tangential trace} is imposed on $\Gamma_{\mathrm{nor}}$, not on $\Gamma_{\mathrm{tan}}$.

Define $d_{\Gamma_{\mathrm{nor}}}^j$ by closing the exterior derivative in the graph norm
on smooth forms that vanish near $\overline{\Gamma_{\mathrm{nor}}}$. Its domain imposes the
essential trace $\iota^*u=0$ on $\Gamma_{\mathrm{nor}}$. Set
\[
\delta_{\Gamma_{\mathrm{tan}}}^j=(d_{\Gamma_{\mathrm{nor}}}^{j-1})^*,\qquad
V_{q,\mathrm{mix}}=D(d_{\Gamma_{\mathrm{nor}}}^q)\cap D(\delta_{\Gamma_{\mathrm{tan}}}^q).
\]
The adjoint domain imposes the zero normal trace
$\iota^*(\star u)=0$ on $\Gamma_{\mathrm{tan}}$ in its weak sense. The densely defined
closed nonnegative form
\begin{equation}\label{eq:mixed-form}
\mathfrak a_{q,\mathrm{mix}}(u,v)
 =\ip{d_{\Gamma_{\mathrm{nor}}}^qu}{d_{\Gamma_{\mathrm{nor}}}^qv}
  +\ip{\delta_{\Gamma_{\mathrm{tan}}}^qu}{\delta_{\Gamma_{\mathrm{tan}}}^qv},
\qquad u,v\in V_{q,\mathrm{mix}},
\end{equation}
defines its self-adjoint operator unambiguously by
\begin{equation}\label{eq:mixed-operator-domain}
\begin{aligned}
D(\Delta_{q,\mathrm{mix}})
 &=\{u\in V_{q,\mathrm{mix}}:\text{ there is }f\in L^2\Omega^q(M)
       \text{ with }\mathfrak a_{q,\mathrm{mix}}(u,v)=\ip{f}{v}
       \text{ for all }v\in V_{q,\mathrm{mix}}\},\\
\Delta_{q,\mathrm{mix}}u&=f.
\end{aligned}
\end{equation}
For smooth operator-domain forms away from the interface, integration by parts
recovers the full local conditions
\begin{align}
\iota^*(\star u)=0,\quad \iota^*(\star du)=0
   &\qquad\text{on }\Gamma_{\mathrm{tan}},\label{eq:mixed-absolute}\\
\iota^*u=0,\quad \iota^*(\delta u)=0
   &\qquad\text{on }\Gamma_{\mathrm{nor}}.\label{eq:mixed-relative}
\end{align}
These conditions cancel the boundary energy away from the interface.
The form construction specifies the realization at the interface, where
singularities may prevent global $H^2$ regularity.

Mixed compactness and closed-range results give a finite-dimensional harmonic
space and a positive spectral gap on its orthogonal complement for admissible
partitions; see \cite[Theorems 4.8--4.9]{PaulySchomburg2022} for bounded
strong Lipschitz pairs. 
For mixed-boundary Hodge theory in the Riemannian setting, see \cite{GoldshteinMitreaMitrea2011}.
In the smooth manifold setting, the corresponding
mixed de Rham complex represents the pair
$(M,\overline{\Gamma_{\mathrm{nor}}})$, so
\begin{equation}\label{eq:mixed-cohomology}
\Ker\Delta_{q,\mathrm{mix}}\simeq
 H_{\mathrm{dR}}^q(M,\overline{\Gamma_{\mathrm{nor}}}).
\end{equation}
The relative boundary piece, as well as the manifold, therefore determines
the kernel; compare \cite[Section 7.1]{Licht2019} and the compatible discrete
construction in \cite[Section 5.3]{ZhaoEtAl2019}. The arguments of
Theorems~\ref{thm:det-invariance}--\ref{thm:decay} and
Proposition~\ref{prop:det-energy} then apply with the mixed operator and its
own harmonic projector and spectral gap under these admissibility assumptions.

\paragraph{Nonnegative boundary dissipation.}
Boundary cancellation is not necessary for a stable diffusion. One controlled
Robin-type extension starts from a standard closed form
$\mathfrak a_{q,B}$ on $V_{q,B}$, a trace map
$T:V_{q,B}\to L^2(\partial M;E)$ continuous in the form norm, and a bounded
nonnegative self-adjoint operator $R$ on this boundary trace space. Here $E$
is the appropriate boundary form bundle and $T$ selects an available trace.
Define
\begin{equation}\label{eq:robin-form}
\mathfrak a_R(u,v)=\mathfrak a_{q,B}(u,v)
 +\ip{R^{1/2}Tu}{R^{1/2}Tv}_{\partial M},\qquad u,v\in V_{q,B}.
\end{equation}
The new form norm is equivalent to the original form norm: the additional
term is nonnegative and bounded in that norm. Thus this densely defined form
is closed and nonnegative, and defines a self-adjoint operator $\Delta_R$.
Its strong heat solutions satisfy
\begin{equation}\label{eq:robin-energy}
\frac12\frac{\dt}{\dt t}\norm{\omega(t)}^2
 =-\nu\left(\norm{d_B^q\omega(t)}^2
           +\norm{\delta_B^q\omega(t)}^2
           +\norm{R^{1/2}T\omega(t)}_{\partial M}^2\right).
\end{equation}
The additional term dissipates energy. However, its zero-energy space is
\[
\Ker\Delta_R
 =\{h\in\Harm_B^q:R^{1/2}Th=0\},
\]
which can be smaller than $\Harm_B^q$. The Robin-type flow preserves
this modified kernel, not necessarily the original cohomology. The
construction depends on the stated trace-continuity and closed-form
assumptions.

\section{Stochastic Hodge--Laplacian diffusion}\label{sec:stoch}
\subsection{Well-posedness and harmonic noise}
Let $U$ be a separable Hilbert space and let $(\Omega,\mathcal F,(\mathcal F_t)_{t\ge0},\mathbb P)$ satisfy the usual conditions. 
Let $W$ be a $U$-cylindrical Wiener process with respect to $(\mathcal F_t)$.
Write $\mathcal L_2(U,H)$ for the Hilbert--Schmidt operators from $U$ to $H$.
Assume $\omega_0\in L^2(\Omega;H)$ is $\mathcal F_0$-measurable and that the given $\mathcal L_2(U,H)$-valued process $G_t$ is predictable and satisfies
\begin{equation}\label{eq:G-HS}
\E\int_0^T\norm{G_s}_{\HS(U,H)}^2\,\dt s<\infty
\qquad\text{for every finite }T.
\end{equation}
The It\^o equation and its mild solution are
\begin{align}
\dt\omega_t&=-\nu\Delta_B\omega_t\,\dt t+G_t\,\dt W_t,
\label{eq:spde}\\
\omega_t&=S_B(t)\omega_0+\int_0^tS_B(t-s)G_s\,\dt W_s.
\label{eq:spde-mild}
\end{align}
Standard variational and semigroup theory \cite{DaPratoZabczyk2014} yields a unique
solution in
$L^2(\Omega;C([0,T];H))\cap L^2(\Omega;L^2(0,T;V_{1,B}))$.
The stochastic differential $\dt\omega_t$ is distinct from the spatial derivative
$d_B\omega_t$.

\begin{theorem}[Characterization of harmonic preservation]\label{thm:stoch-invariance}
Under \eqref{eq:G-HS},
\begin{equation}\label{eq:harmonic-noise-identity}
\E\norm{P_B(\omega_t-\omega_0)}^2
 =\E\int_0^t\norm{P_BG_s}_{\HS(U,H)}^2\,\dt s.
\end{equation}
In particular, harmonic coordinates are preserved almost surely for all times if and
only if $P_BG_s=0$ for $\dt s\otimes\mathbb P$-almost every $(s,\omega)$.
\end{theorem}
\begin{proof}
Applying the bounded projector to \eqref{eq:spde-mild} and using $P_BS_B(t)=P_B$
gives $P_B(\omega_t-\omega_0)=\int_0^tP_BG_s\,\dt W_s$.
It\^o isometry proves the identity and its necessity. Sufficiency holds simultaneously
at all times by continuity, first intersecting probability-one events over rational times.
\end{proof}
The canonical constrained forcing is $G_t=Q_B\Sigma_t$, where $\Sigma$
is predictable and satisfies \eqref{eq:G-HS}. For an orthonormal
harmonic basis $\{h_j\}_{j=1}^r$, the coefficient
$c_j(t)=\ip{\omega_t}{h_j}$ has quadratic variation
$\int_0^t\norm{G_s^*h_j}_U^2\,\dt s$. The Hodge drift provides no restoring force
in these directions. 
Although $\E c_j(t)=\E c_j(0)$, individual trajectories need not preserve the harmonic coefficient when $G_t^*h_j\ne0$.

\subsection{Energy and stationary fluctuations}
Set $Y_t=Q_B\omega_t$ and $\mathcal I(t)=\E\norm{Q_BG_t}_{\HS(U,H)}^2$.
The variational It\^o identity gives
\begin{equation}\label{eq:ito-energy-rigorous}
\E\norm{Y_t}^2+2\nu\E\int_0^t\norm{\Delta_B^{1/2}Y_s}^2\,\dt s
 =\E\norm{Y_0}^2+\int_0^t\mathcal I(s)\,\dt s.
\end{equation}
Combining it with \eqref{eq:poincare} yields
\begin{equation}\label{eq:stochastic-energy-convolution}
\E\norm{Y_t}^2\le e^{-2\nu\lambda_{B,+}t}\E\norm{Y_0}^2
 +\int_0^te^{-2\nu\lambda_{B,+}(t-s)}\mathcal I(s)\,\dt s.
\end{equation}
If $\mathcal I(t)\le \mathcal I_*$, then
$\limsup_{t\to\infty}\E\norm{Y_t}^2\le \mathcal I_*/(2\nu\lambda_{B,+})$.
The deterministic relaxation thus coexists with nonzero stochastic fluctuations.

For the invariant-law discussion, assume $G_t\equiv G\in\mathcal L_2(U,H)$ is deterministic and $P_BG=0$. 
The Ornstein--Uhlenbeck process on $Q_BH$ has invariant
law $\mathcal N(0,C_\infty)$ with
\begin{equation}\label{eq:invariant-covariance}
C_\infty=\int_0^\infty S_B(s)GG^*S_B(s)\,\dt s,
\qquad \nu(\Delta_BC_\infty+C_\infty\Delta_B)=GG^*.
\end{equation}
The operator equation is understood weakly on $D(\Delta_B)\cap Q_BH$. Its covariance is trace class
because $\tr C_\infty\le\norm{G}_{\HS}^2/(2\nu\lambda_{B,+})$.
Exponential contraction makes this invariant probability law unique on $Q_BH$ even
if $G$ is degenerate; nondegeneracy is not required for uniqueness. Synchronously
coupled solutions differ by $S_B(t)(Y_0-\widetilde Y_0)$, and the mild solution
converges in law to the same Gaussian from every fixed initial state.

On the full space, the harmonic component is conserved. For each
$h_0\in\Harm_B$, the translate of $\mathcal N(0,C_\infty)$ by $h_0$
is invariant on $h_0+Q_BH$. Mixtures over $h_0$ are invariant on $H$.
If $\dim\Harm_B>0$, there is therefore no unique invariant
probability law across all harmonic coordinates. If a constant additive $G$ has $P_BG\ne0$, its harmonic
projection is a nontrivial Brownian motion and has no invariant probability law.

For a diagonal-noise example, choose an orthonormal basis $\{u_j\}$ of $U$
indexed by the positive orthonormal eigenbasis
$\Delta_Be_j=\lambda_je_j$, $\lambda_j>0$, and set $Gu_j=\sigma_je_j$.
Then the modal coefficient $a_j(t)=\ip{Q_B\omega_t}{e_j}$ satisfies
\begin{equation}\label{eq:mode-variance}
\E|a_j(t)|^2=e^{-2\nu\lambda_jt}\E|a_j(0)|^2
 +\frac{\sigma_j^2}{2\nu\lambda_j}(1-e^{-2\nu\lambda_jt}).
\end{equation}
This is a second-moment identity; the corresponding variance formula replaces both second moments by variances. 
It provides an exact reference for numerical experiments. 
Our baseline assumes $\sum_j\sigma_j^2<\infty$. 
Rougher noise requires a separate stochastic-convolution
argument; $\sum_j\sigma_j^2/\lambda_j<\infty$ is the stationary finite-energy condition
in the diagonal case. Projecting out a finite-dimensional harmonic space alone does
not make spatial white noise Hilbert--Schmidt.

\section{Topology-preserving random perturbations}\label{sec:random}
Section~\ref{sec:stoch} gives a constraint on continuous-time noise. The same
geometry also characterizes random perturbations of individual data samples and
clarifies which perturbations preserve a cohomology class.

\subsection{Projected, exact, and coexact perturbations}
For a square-integrable $H$-valued random variable $\zeta$, the perturbation
\begin{equation}\label{eq:random-projection}
\widetilde\omega=\omega+Q_B\zeta
\end{equation}
preserves $P_B\omega$ for every realization. It can alter both exact and coexact
components. More generally, for degree $q$, a potential construction is
\begin{equation}\label{eq:potential-perturbation}
\widetilde\omega=\omega+d_B^{q-1}\psi+\delta_B^{q+1}\eta,
\end{equation}
where the random potentials are measurable and square-integrable in
the respective graph norms. Orthogonality
in \eqref{eq:hodge-decomp} makes both added terms harmonic-free. Similarly,
continuous forcing $G_t=d_B^{q-1}K_t+\delta_B^{q+1}J_t$ is harmonic-free whenever
the potential-valued operators are square-integrable Hilbert--Schmidt in these
graph norms. These constructions expose the types of variation introduced by the
noise, whereas an arbitrary orthogonal projection need not distinguish them.

If the original form is closed, exact perturbations preserve its cohomology class.
Coexact perturbations generally destroy closedness. For example, in dimension $m\ge2$, choose a smooth
compactly supported $2$-form $\eta$ in the interior with $d\delta\eta\ne0$; then
$\delta\eta$ is orthogonal to the harmonic space for either boundary realization,
but adding it to a closed $1$-form produces a nonclosed form. Thus zero harmonic
leakage is insufficient to justify a cohomology claim.

\subsection{Cohomology-preserving stochastic diffusion}
We state the result for any degree $1\le q\le m$. Equip $D(d_B^{q-1})$ with the Hilbert graph norm
\[
\|\alpha\|_{D(d_B^{q-1})}^2
=\|\alpha\|_{L^2}^2+\|d_B^{q-1}\alpha\|_{L^2}^2.
\]
Assume
\begin{equation}\label{eq:K-assumption}
\E\int_0^T\norm{K_s}_{\HS(U,D(d_B^{q-1}))}^2\,\dt s<\infty,
\end{equation}
for every finite $T$, where $K$ is predictable. This makes $G_s=d_B^{q-1}K_s$ a valid noise operator
under \eqref{eq:G-HS} in degree $q$.

\begin{theorem}[Preservation of absolute or relative cohomology]
\label{thm:cohomology-preservation}
Let $\omega_0\in L^2(\Omega;Z_B^q)$ be $\mathcal F_0$-measurable and suppose
\eqref{eq:K-assumption} holds. The mild solution of
\begin{equation}\label{eq:exact-noise}
\dt\omega_t=-\nu\Delta_{q,B}\omega_t\,\dt t+d_B^{q-1}K_t\,\dt W_t
\end{equation}
satisfies, on one probability-one event for every $t\ge0$,
\begin{equation}\label{eq:class-preservation}
\omega_t\in P_B^q\omega_0+E_B^q,\qquad
d_B^q\omega_t=0,\qquad P_B^q\omega_t=P_B^q\omega_0,\qquad
[\omega_t]=[\omega_0].
\end{equation}
The class belongs to $H_{\mathrm{dR}}^q(M)$ for $B=\mathrm{tan}$ and to
$H_{\mathrm{dR}}^q(M,\partial M)$ for $B=\mathrm{nor}$. For nonsmooth forms it is
interpreted through the corresponding $L^2$ Hilbert complex.
\end{theorem}
\begin{proof}
Write $\omega_0=h_0+e_0$ using \eqref{eq:closed-decomp}. The mild formula gives
\[
\omega_t=h_0+S_{q,B}(t)e_0+
 \int_0^tS_{q,B}(t-s)d_B^{q-1}K_s\,\dt W_s.
\]
The exact space $E_B^q$ is closed and invariant under the heat semigroup. Both
non-harmonic terms therefore belong to $E_B^q$; for the stochastic integral this
follows by approximation with $E_B^q$-valued finite sums. Continuous paths and
closedness of $E_B^q$ make the assertion simultaneous for all times. Consequently
$\omega_t-\omega_0\in E_B^q$, which proves all claims.
\end{proof}

The proof uses closed-subspace invariance rather than differentiating
the stochastic integral. Closed range supplies a bounded minimal-norm
right inverse $R_B:E_B^q\to D(d_B^{q-1})$, giving the potential
$R_B(\omega_t-\omega_0)$ for the exact difference.

Relative potentials have zero tangential trace weakly (Dirichlet data
for scalars); absolute potentials belong to the maximal domain without
an essential trace. These are the domains specified in
Section~\ref{sec:hodge}, not a common vector-field boundary condition
imposed on every degree.

Boundary conditions hold in the mild/variational sense; classical
second-order traces are not asserted for every sample path. A sufficient
stronger setting is $\omega_0\in L^2(\Omega;D(\Delta_{q,B}))$ and
\eqref{eq:G-HS} with $H$ replaced by $D(\Delta_{q,B})$, equipped with
its graph norm. The corresponding solution has continuous paths in that
graph space.

\subsection{Constructive examples and observables}
For a concrete degree-one construction, let $\{u_\ell\}$ be an orthonormal basis of $U$, and let $\{\phi_\ell\}$ be orthonormal scalar eigenfunctions satisfying $\Delta_{0,B}\phi_\ell=\mu_\ell\phi_\ell$ with $\mu_\ell>0$.
Use Neumann conditions for $B=\mathrm{tan}$ and Dirichlet conditions for $B=\mathrm{nor}$.
Then
\begin{equation}\label{eq:scalar-exact-noise}
\theta_\ell=\mu_\ell^{-1/2}d\phi_\ell,\qquad
Gu_\ell=\sigma_\ell\theta_\ell,\qquad
Ku_\ell=\sigma_\ell\mu_\ell^{-1/2}\phi_\ell
\end{equation}
gives orthonormal exact $1$-form modes with
$\Delta_{1,B}\theta_\ell=\mu_\ell\theta_\ell$. They satisfy both second-order
boundary conditions. For the normal realization, for example,
$\iota^*d\phi_\ell=d\iota^*\phi_\ell=0$ and
$\iota^*\delta d\phi_\ell=\mu_\ell\iota^*\phi_\ell=0$.
The Neumann condition gives the two tangential traces analogously. If
$\sum_\ell\sigma_\ell^2<\infty$, the scalar spectral gap verifies
\eqref{eq:K-assumption}; finite-rank choices give simple smooth benchmarks.

For sufficiently smooth closed forms, periods provide basis-independent cohomology
observables: $\int_\gamma\omega_t$ is unchanged for every absolute cycle $\gamma$,
or relative cycle with boundary in $\partial M$ in the relative case. Stokes' theorem
applied to the exact difference proves this. These cycle integrals are not automatically
defined for arbitrary $L^2$ forms; the Hilbert-complex statement
\eqref{eq:class-preservation} remains the primary analytical assertion. On discrete
cochains the cycle pairings are directly available.

For a semilinear It\^o extension with additional drift $F_t$,
harmonic preservation requires $P_B^qF_t=0$ and $P_B^qG_t=0$ along the solution.
For a prescribed class $h_0+E_B^q$, with $h_0\in\Harm_B^q$,
sufficient invariance conditions are $F_t\in E_B^q$ and
$G_t(U)\subseteq E_B^q$ on that class, together with assumptions
ensuring well-posedness. 
The learned reverse model in Section~\ref{sec:reverse} enforces the corresponding finite-dimensional restriction through its coordinates.

\section{Cellular sheaf discretization and topology-preserving dynamics}\label{sec:discrete}
\subsection{The constant cellular sheaf and the boundary complex}
A matrix kernel can be invariant without representing the intended cohomology.
Let $K_h$ be a boundary-conforming finite regular cell complex modeling $M$, with
boundary subcomplex $\partial K_h$, and let $\constsheaf_{K_h}$ be its constant
cellular sheaf: every cell carries $\R$ and every restriction map
is the identity. Its cochains are ordinary scalar cellular cochains. Take
$C_{h,\mathrm{tan}}^j=C^j(K_h;\constsheaf_{K_h})$ for the absolute realization and
$C_{h,\mathrm{nor}}^j=C^j(K_h,\partial K_h;\constsheaf_{K_h})$ for the relative
realization, with relative cochains vanishing on the boundary subcomplex.
The real coordinates below are column vectors.

Write $d_{h,B}^j:C_{h,B}^j\to C_{h,B}^{j+1}$ for the abstract differential
and $D_B^j$ for its coordinate matrix. The superscript is the source degree,
not a power. Thus
\begin{equation}\label{eq:discrete-complex}
d_{h,B}^{j+1}d_{h,B}^j=0,\qquad D_B^{j+1}D_B^j=0.
\end{equation}
The matrices $D_{\mathrm{tan}}^j,D_{\mathrm{nor}}^j$ are the constant-sheaf coboundaries with the
selected boundary treatment. 

\begin{proposition}[Constant-sheaf identification]\label{prop:constant-sheaf-identification}
For $B=\mathrm{tan}$, the cohomology of \eqref{eq:discrete-complex} is
$H^q(K_h;\constsheaf)$; for $B=\mathrm{nor}$, it is
$H^q(K_h,\partial K_h;\constsheaf)$. With positive definite mass inner
products, the operator whose matrices are specified in
\eqref{eq:discrete-adjoint}--\eqref{eq:discrete-stiffness}
is the Hodge Laplacian of this constant-sheaf or relative cochain complex.
When the masses are orthogonal sums of stalk inner products, this is the
weighted cellular sheaf Laplacian in the sense of the cited sheaf literature.
\end{proposition}
\begin{proof}
Each restriction is the identity, so \eqref{eq:sheaf-coboundary} reduces to
signed cellular incidence. Restricting to cochains vanishing on the boundary
gives the relative subcomplex. The adjoint and Laplacian are determined by
the chosen mass inner products.
\end{proof}

For a Cartesian level-set implementation, the selected normal or tangential support is not automatically a boundary-conforming cell complex. 
In particular, the normal inclusion rule can retain an edge incident to an interior vertex while omitting its exterior endpoint \cite[Section~3.2]{PersistentHodge2024}.
The support-restricted coboundary operators in \cite[Section~3.2]{PersistentHodge2024} satisfy the complex identity under the support-compatibility conditions established there. 
For the stochastic annulus experiment, we also verified $D^1D^0=0$ exactly for the assembled matrices (Section~\ref{sec:num-perturbations}).

Applying Proposition~\ref{prop:constant-sheaf-identification} to such arrays requires an explicit complex or relative-pair interpretation of the retained degrees of freedom, together with the matching incidence and mass matrices.
The semidiscrete algebra below still applies to any verified finite Hilbert complex, but its cohomology must not be identified with that of $M$ solely from nilpotency or an observed kernel dimension.

Let $M_{j,B}$ be symmetric positive definite mass matrices, with the mesh
index suppressed. In diagonal DEC these are the discrete Hodge-star matrices of
\cite{ZhaoEtAl2019,LiuEtAl2026}; $M$ also accommodates Galerkin masses. The orthogonal-direct-sum convention of
\cite{WeiWei2025PSL,HansenGhrist2019} gives block-diagonal stalk masses.
A general non-block-diagonal Galerkin mass is an additional global inner
product on the same cochain complex, not merely a change of individual stalk
weights. The algebraic Hodge formulas below allow either choice.
Define the downward adjoint matrix $\delta_{j,h,B}$ and the
raw-coordinate Laplacian matrix $A_{q,h,B}$ by
\begin{equation}\label{eq:discrete-adjoint}
\begin{aligned}
\delta_{j,h,B}&=M_{j-1,B}^{-1}(D_B^{j-1})^{\top} M_{j,B},\\
A_{q,h,B}&=D_B^{q-1}\delta_{q,h,B}+\delta_{q+1,h,B}D_B^q.
\end{aligned}
\end{equation}
The corresponding symmetric stiffness matrix is denoted $\mathsf K_{q,B}$:
\begin{equation}\label{eq:discrete-stiffness}
\begin{aligned}
\mathsf K_{q,B}&=M_{q,B}A_{q,h,B}\\
 &=M_{q,B}D_B^{q-1}M_{q-1,B}^{-1}(D_B^{q-1})^{\top} M_{q,B}
   +(D_B^q)^{\top} M_{q+1,B}D_B^q.
\end{aligned}
\end{equation}
For diagonal DEC masses, this is the symmetric stiffness matrix used in
\cite[Section~3.2, equations~(3.9)--(3.10)]{PersistentHodge2024}.
We reserve $L$ for the Laplacian matrix in orthonormal coordinates.
Specifically, set
\begin{equation}\label{eq:orthonormal-laplacian}
\begin{aligned}
\widetilde D_B^j&=M_{j+1,B}^{1/2}D_B^jM_{j,B}^{-1/2},\\
L_{q,h,B}&=M_{q,B}^{-1/2}\mathsf K_{q,B}M_{q,B}^{-1/2}\\
 &=\widetilde D_B^{q-1}(\widetilde D_B^{q-1})^{\top}
    +(\widetilde D_B^q)^{\top}\widetilde D_B^q.
\end{aligned}
\end{equation}
The square roots are the symmetric positive definite ones. Thus $A_{q,h,B}$
is self-adjoint in the $M_{q,B}$ inner product, whereas $L_{q,h,B}$ is symmetric
in the Euclidean inner product. They represent the same operator in different
coordinates and have the same spectrum; the unweighted spectrum of
$\mathsf K_{q,B}$ is not in general that spectrum. Equivalently,
\[
 \mathsf K_{q,B}v=\lambda M_{q,B}v
 \quad\Longleftrightarrow\quad
 L_{q,h,B}(M_{q,B}^{1/2}v)=\lambda M_{q,B}^{1/2}v.
\]
When all masses are the identity, $\mathsf K_{q,B}=A_{q,h,B}=L_{q,h,B}$,
and the ordinary transpose formulas of \cite{WangWei2026MPSL,CottrellEtAl2026AD}
apply without a coordinate change. In arbitrary coordinates an adjoint is
\emph{not} obtained by merely transposing $D_B^j$; the Gram-matrix dependence
is also explicit in \cite[Proposition~3.1]{WeiWei2025PSL}.
The energy identity is
\begin{equation}\label{eq:discrete-energy}
x^\top\mathsf K_{q,B}x=\norm{\delta_{q,h,B}x}_{M_{q-1,B}}^2
                         +\norm{D_B^q x}_{M_{q+1,B}}^2.
\end{equation}
Finite-dimensional Hodge theory gives
\begin{equation}\label{eq:discrete-closed-split}
\begin{aligned}
\Ker D_B^q&=\Range D_B^{q-1}\oplus_{M_{q,B}}\Ker A_{q,h,B},\\
\dim\Ker A_{q,h,B}&=\dim(\Ker D_B^q/\Range D_B^{q-1}).
\end{aligned}
\end{equation}
These dimensions agree with \eqref{eq:nullity} when the complex represents
the intended manifold or pair. Continuum approximation of the metric and
operator is an additional requirement. DEC supplies relevant constructions
\cite{ZhaoEtAl2019,SuTongWei2024}, while FEEC explains the role of
subcomplexes and bounded commuting projections in stability and convergence
\cite{ArnoldFalkWinther2010}. These algebraic identities alone do not prove a
continuum convergence rate.

Support-selection matrices in
\cite[Section~3.2]{PersistentHodge2024} are distinct from harmonic
projectors. Identities for one complex do not automatically hold for
transfers between different boundary supports.

For a fixed weight, replace all masses by the weighted masses and recompute
adjoints and harmonic projectors \cite{SuTongWei2026}. The incidence complex
and its cohomology are unchanged. In diagonal DEC,
$M_{j,B,f}=\diag(e^{-f_j})M_{j,B}$, where $f_j$ denotes the sampled
weight function on the degree-$j$ cells.
A constant $f=c_f$ scales both $M_{q,B}$ and $\mathsf K_{q,B}$ by $e^{-c_f}$,
leaving $A_{q,h,B}$, $L_{q,h,B}$, and the generalized spectrum unchanged.
Nonconstant weights generally change the harmonic representatives.

\subsection{General weighted cellular sheaves}\label{sec:general-cellular-sheaf}
Let $\cS_h$ be a weighted cellular sheaf on $K_h$. Use its full cochain
complex or its relative complex for a specified boundary subcomplex.
Write $d_{\cS_h}^j$ for its coboundary and $D_{\cS_h}^j$ for its matrix,
assembled from the signed restrictions in \eqref{eq:sheaf-coboundary}.
Let $M_{j,\cS_h}$ be positive definite block mass matrices representing
stalk inner products. Define
\begin{equation}\label{eq:general-sheaf-adjoint}
 \delta_{j,\cS_h}=M_{j-1,\cS_h}^{-1}(D_{\cS_h}^{j-1})^{\top} M_{j,\cS_h}.
\end{equation}
Then
\begin{equation}\label{eq:general-sheaf-laplacian}
\begin{aligned}
 A_{q,\cS_h}&=D_{\cS_h}^{q-1}\delta_{q,\cS_h}
          +\delta_{q+1,\cS_h}D_{\cS_h}^q,\\
 \mathsf K_{q,\cS_h}&=M_{q,\cS_h}A_{q,\cS_h},\qquad
 L_{q,\cS_h}=M_{q,\cS_h}^{-1/2}\mathsf K_{q,\cS_h}M_{q,\cS_h}^{-1/2}.
\end{aligned}
\end{equation}
The last matrix is the orthonormal-coordinate representation of
$\Delta_{q,\cS_h}$. In the full complex, sheaf Hodge theory yields
\begin{equation}\label{eq:general-sheaf-hodge}
\begin{aligned}
\Ker A_{q,\cS_h}&\simeq H^q(K_h;\cS_h),\\
\Ker D_{\cS_h}^q&=\Range D_{\cS_h}^{q-1}
 \oplus_{M_{q,\cS_h}}\Ker A_{q,\cS_h}.
\end{aligned}
\end{equation}
Use relative sheaf cohomology instead for the relative complex.
The coefficient sheaf, restriction maps, and mass matrices are held fixed during each evolution considered below. 
All algebraic preservation statements below use the Hilbert-complex identities, self-adjointness, and the harmonic projector, and remain valid with these sheaf-valued matrices. 
The constant-sheaf case is retained in the main notation to keep the continuum comparison explicit.

\subsection{Semidiscrete preservation}
Fix degree $q$ and boundary type $B$. In the evolution equations write
$M_h=M_{q,B}$, $\mathsf K_h=\mathsf K_{q,B}$, and $A_h=M_h^{-1}\mathsf K_h$,
and use $\norm{x}_{M_h}^2=x^\top M_hx$. Suppress the fixed $B$ in
$D^j=D_B^j$, $M_j=M_{j,B}$, and $\delta_{j,h}=\delta_{j,h,B}$.
The orthonormal-coordinate Laplacian is
$L_h=L_{q,h,B}=M_h^{-1/2}\mathsf K_hM_h^{-1/2}$.
The potential-noise matrix $\mathcal K_h$ is distinct from the base complex $K_h$
and the stiffness $\mathsf K_h$. Let
\begin{equation}\label{eq:projector-discrete}
\mathsf K_hH_h=0,\quad H_h^\top M_hH_h=I,\qquad
P_h=H_hH_h^\top M_h,\quad Q_h=I-P_h.
\end{equation}
Here the columns of $H_h$ span the entire kernel of the assembled matrix;
computed approximate subspaces are distinguished below. Write
$N_h=\dim C_{h,B}^q$. For a fixed deterministic matrix
$B_h\in\R^{N_h\times r_{\mathrm W,h}}$ and standard $r_{\mathrm W,h}$-dimensional Brownian motion, consider
\begin{equation}\label{eq:semi-discrete}
\dt X_t=-\nu A_hX_t\,\dt t+B_h\,\dt W_t.
\end{equation}
Use the filtration of Section~\ref{sec:stoch}, with $W$ a Brownian
motion relative to it and $X_0\in L^2(\Omega;\R^{N_h})$
$\mathcal F_0$-measurable. All matrices are fixed. In mass-orthonormal
coordinates, set
$\widetilde X_t=M_h^{1/2}X_t$ and $\widetilde B_h=M_h^{1/2}B_h$. Then the
same dynamics and harmonic projector are
\begin{equation}\label{eq:orthonormal-sde}
\begin{aligned}
 \dt\widetilde X_t&=-\nu L_h\widetilde X_t\,\dt t
                       +\widetilde B_h\,\dt W_t,\\
 \widetilde P_h&=M_h^{1/2}P_hM_h^{-1/2}
               =\widetilde H_h\widetilde H_h^\top,
 \qquad \widetilde H_h=M_h^{1/2}H_h.
\end{aligned}
\end{equation}
This coordinate change preserves the constraint:
$P_hB_h=0$ if and only if $\widetilde P_h\widetilde B_h=0$.

\begin{proposition}[Discrete harmonic and class preservation]
\label{prop:discrete-kernel-invariance}
If $P_hB_h=0$, then $P_hX_t=P_hX_0$ almost surely for all times.
If additionally $D^qX_0=0$ and $B_h=D^{q-1}\mathcal K_h$, then
\begin{equation}\label{eq:discrete-class}
X_t\in P_hX_0+\Range D^{q-1},\qquad
D^qX_t=0,\qquad [X_t]=[X_0]\text{ in }\Ker D^{q}/\Range D^{q-1}.
\end{equation}
\end{proposition}
\begin{proof}
Symmetry gives $P_hA_h=H_hH_h^\top \mathsf K_h=0$, proving the first assertion by
projection. For the second, $A_h$ maps the exact range to itself by
\eqref{eq:discrete-complex}, and \eqref{eq:discrete-closed-split} decomposes the
initial state into harmonic and exact parts. The matrix-exponential solution and
its stochastic convolution therefore stay in the displayed affine space.
\end{proof}

\begin{corollary}[Sheaf-valued preservation]\label{cor:sheaf-valued-preservation}
Let $\cS_h$ be any weighted cellular sheaf with Laplacian
$A_{q,\cS_h}$ and harmonic projector $P_{q,\cS_h}$. The semidiscrete sheaf flow
\[
\dt X_t=-\nu A_{q,\cS_h}X_t\,\dt t+B_{\cS_h}\,\dt W_t
\]
preserves the harmonic sheaf-cohomology coordinates whenever
$P_{q,\cS_h}B_{\cS_h}=0$. If the initial cochain is closed and
$B_{\cS_h}=D_{\cS_h}^{q-1}\mathcal K_h$, then its class in $H^q(K_h;\cS_h)$ is preserved.
For a relative complex, replace this group by the corresponding relative sheaf
cohomology throughout.
\end{corollary}

Noise amplitudes must be specified with their coordinates or covariance
factor. Independent unit-variance raw coefficients do not, by themselves,
approximate a fixed continuum $L^2$ forcing under mesh refinement.
Projected and exact constructions also need not inject the same covariance.

\subsection{Time stepping: invariance, stability, and covariance}
Set $t_n=n\tau$ and $\Delta W_n=W_{t_{n+1}}-W_{t_n}$.
Euler--Maruyama \cite{Higham2001} uses these independent increments,
with $\Delta W_n\sim\mathcal N(0,\tau I)$, in the update
\begin{equation}\label{eq:em-step}
X^{n+1}=R_\tau X^n+B_h\Delta W_n,\qquad R_\tau=I-\nu\tau A_h.
\end{equation}
In implementation, solve with $M_h$ rather than explicitly forming its inverse.
Since $P_hR_\tau=P_h$, the update preserves harmonic coefficients for every time
step in exact arithmetic when $P_hB_h=0$. Exact noise and closed initial data also
preserve \eqref{eq:discrete-class} at each step. These are algebraic properties,
independent of the step's stability.

If the positive spectrum is empty, a harmonic-preserving update is constant and
needs no stability restriction. Otherwise, let $\lambda_{h,\max}$ be the largest positive generalized
eigenvalue of $\mathsf K_hv=\lambda M_hv$. With
\begin{equation}\label{eq:em-cfl}
0<\nu\tau\lambda_{h,\max}<2,\qquad
\rho_\tau=\max_{\lambda_j>0}|1-\nu\tau\lambda_j|<1,
\end{equation}
the positive modes contract. For $Y^n=Q_hX^n$ and
$b_{0,h}^2=\tr[(Q_hB_h)^{\top} M_h(Q_hB_h)]$, conditional expectation gives
\begin{align}
\E[\norm{Y^{n+1}}_{M_h}^2\mid\mathcal F_{t_n}]
 &=\norm{R_\tau Y^n}_{M_h}^2+\tau b_{0,h}^2,\label{eq:em-conditional-energy}\\
\E\norm{Y^n}_{M_h}^2
 &\le\rho_\tau^{2n}\E\norm{Y^0}_{M_h}^2
       +\tau b_{0,h}^2\frac{1-\rho_\tau^{2n}}{1-\rho_\tau^2}.
       \label{eq:em-mean-square-bound}
\end{align}
At equality in the threshold, a highest mode can alternate without damping, and
noise in that mode accumulates variance. Above it, some positive modes grow even
though harmonic coefficients remain fixed.

For diagonal modal noise of amplitude $\sigma_j$, the stationary explicit variance is
\begin{equation}\label{eq:em-stationary-variance}
v_{j,\tau}=\frac{\tau\sigma_j^2}{1-(1-\nu\tau\lambda_j)^2}
 =\frac{\sigma_j^2}{\nu\lambda_j(2-\nu\tau\lambda_j)}.
\end{equation}
It exceeds the continuous variance $\sigma_j^2/(2\nu\lambda_j)$ whenever
$\sigma_j\ne0$. The exact-in-law modal reference update is
\begin{equation}\label{eq:exact-ou-reference-step}
a_j^{n+1}=e^{-\nu\lambda_j\tau}a_j^n+
 \sigma_j\sqrt{\frac{1-e^{-2\nu\lambda_j\tau}}{2\nu\lambda_j}}\,\xi_{j,n},
\qquad \xi_{j,n}\sim\mathcal N(0,1).
\end{equation}
For the diagonal-noise reference, the variables $\xi_{j,n}$ are independent standard normals.
The harmonic component remains unchanged when $P_hB_h=0$; otherwise its
increment is $P_hB_h\Delta W_n$. 
Independently simulated reference paths validate distributions, not strong errors against another scheme.

A drift-implicit alternative is
\begin{equation}
 (I+\nu\tau A_h)X^{n+1}=X^n+B_h\Delta W_n.
 \label{eq:implicit-em-step}
\end{equation}
Here the resolvent $(I+\nu\tau A_h)^{-1}$ acts on both the
previous state and the noise increment. For diagonal noise of amplitude
$\sigma_j$ in a mass-orthonormal eigenbasis, a mode with
$\lambda_j>0$ satisfies
\[
 a_j^{n+1}
 = \frac1{1+\nu\tau\lambda_j}a_j^n + \frac{\sigma_j\sqrt{\tau}}{1+\nu\tau\lambda_j}\xi_{j,n}  ,   
 \qquad \xi_{j,n}\sim\mathcal N(0,1),
\]
where the innovations are independent across time and independent of
the initial data. Its stationary variance is therefore
\[
 v_{j,\tau}^{\mathrm{imp}}
 =
 \frac{\sigma_j^2\tau/(1+\nu\tau\lambda_j)^2}
      {1-(1+\nu\tau\lambda_j)^{-2}}
 =
 \frac{\sigma_j^2}{\nu\lambda_j(2+\nu\tau\lambda_j)}.
\]
Adding noise after the implicit solve, without the resolvent factor,
would define a different scheme.

The positive-mode factors $(1+\nu\tau\lambda_j)^{-1}$ lie in $(0,1)$
for every $\tau>0$, so there is no explicit stability threshold.
The scheme preserves the same harmonic and exact-class constraints
under the conditions above, but its stationary variance is below the
continuous value $\sigma_j^2/(2\nu\lambda_j)$ when $\sigma_j\ne0$.
Structure preservation, time stability, and statistical accuracy
must therefore be tested separately.

\subsection{Leakage from approximate harmonic subspaces}\label{sec:topological-leakage}
Let $\widehat P_h$ be a fixed deterministic $M_h$-orthogonal projector.
It need not equal $P_h$. Define
\begin{equation}\label{eq:leakage-defects}
\mathcal L(t)=\norm{\widehat P_h(X_t-X_0)}_{M_h},\qquad
\varepsilon_L=\norm{\widehat P_hA_h}_{M_h\to M_h},\quad
\varepsilon_B=\norm{\widehat P_hB_h}_{\mathrm{HS},M_h},
\end{equation}
where $\norm{B}_{\mathrm{HS},M_h}^2=\tr(B^\top M_hB)$.

\begin{proposition}[Finite-horizon leakage bound]\label{prop:topological-leakage}
For \eqref{eq:semi-discrete},
\begin{equation}\label{eq:leakage-integral-bound}
\E\mathcal L(t)^2\le
 2\nu^2\varepsilon_L^2t\int_0^t\E\norm{X_s}_{M_h}^2\,\dt s
 +2\varepsilon_B^2t.
\end{equation}
If $\sup_{s\le T}\E\norm{X_s}_{M_h}^2\le K_T$, then for $t\le T$,
\begin{equation}\label{eq:leakage-finite-horizon-bound}
\E\mathcal L(t)^2\le2\nu^2K_T\varepsilon_L^2t^2+2\varepsilon_B^2t.
\end{equation}
One may take $K_T=\E\norm{X_0}_{M_h}^2+T\norm{B_h}_{\mathrm{HS},M_h}^2$.
\end{proposition}
\begin{proof}
Project the integral SDE:
$\widehat P_h(X_t-X_0)=-\nu\int_0^t\widehat P_hA_hX_s\,\dt s
 +\int_0^t\widehat P_hB_h\,\dt W_s$.
Apply $\norm{u+v}^2\le2\norm{u}^2+2\norm{v}^2$, Cauchy--Schwarz in time,
and the mass-weighted It\^o isometry. No independence of the drift and stochastic
integrals is assumed. The energy identity for \eqref{eq:semi-discrete} bounds
$\E\norm{X_s}_{M_h}^2$ by its initial value plus
$s\norm{B_h}_{\mathrm{HS},M_h}^2$.
\end{proof}

For a mass-orthonormal candidate basis $\widehat H_h$ with
$\widehat P_h=\widehat H_h\widehat H_h^\top M_h$, the residuals are computable as
\begin{equation}\label{eq:computable-leakage-residuals}
\varepsilon_L=\norm{M_h^{-1/2}\mathsf K_h\widehat H_h}_2,\qquad
\varepsilon_B=\norm{\widehat H_h^\top M_hB_h}_F.
\end{equation}
Here $\|\cdot\|_2$ denotes the matrix spectral norm, and $\|\cdot\|_F$ denotes the Frobenius norm.
For the exact projector the sharper identity is
\begin{equation}\label{eq:exact-kernel-noise-leakage}
\E\norm{P_h(X_t-X_0)}_{M_h}^2=t\norm{P_hB_h}_{\mathrm{HS},M_h}^2.
\end{equation}
Small computed leakage does not certify the subspace's dimension or agreement with
the true harmonic space. The projector discrepancy enters through
\begin{equation}\label{eq:true-projector-leakage-bound}
\begin{aligned}
\E\norm{P_h(X_t-X_0)}_{M_h}^2
 &\le2\E\mathcal L(t)^2\\
 &\quad+4\norm{P_h-\widehat P_h}_{M_h\to M_h}^2(\E\norm{X_t}_{M_h}^2+\E\norm{X_0}_{M_h}^2).
\end{aligned}
\end{equation}
The subspace comparison and topology check are separate from preservation of
$\widehat P_h$.

Summing \eqref{eq:em-step} gives the corresponding time-discrete estimate
\begin{equation}\label{eq:em-leakage-bound}
\E\norm{\widehat P_h(X^n-X^0)}_{M_h}^2
 \le2\nu^2\varepsilon_L^2t_n\tau\sum_{j=0}^{n-1}\E\norm{X^j}_{M_h}^2
       +2\varepsilon_B^2t_n.
\end{equation}
Under $0<\nu\tau\lambda_{h,\max}\le2$, the full deterministic step is a
contraction, so \eqref{eq:leakage-finite-horizon-bound} holds at $t_n\le T$ with
the same $K_T$. 

For the implicit scheme \eqref{eq:implicit-em-step}, the same argument gives \eqref{eq:em-leakage-bound} with the sum over $j=0,\ldots,n-1$ replaced by the sum over $j=1,\ldots,n$.
Its resolvent is contractive for every $\tau>0$, so \eqref{eq:leakage-finite-horizon-bound} holds at $t_n\le T$
with the same $K_T$, without a step-size restriction.

These bounds describe the written schemes in exact arithmetic;
solver and rounding residuals must be added if they are to be bounded.

\section{Reverse diffusion, score models, and generation}\label{sec:reverse}
The forward diffusion supplies a tractable conditional noising process. Its reverse
must be defined on the space actually explored by that process. 
At a fixed cellular discretization, the harmonic coordinates specify a preserved projection. 
For closed cochains, that projection represents their ordinary or relative cohomology class. 
We formulate the generative construction at a fixed finite discretization or finite spectral truncation.
A general infinite-dimensional score or time-reversal
theorem for arbitrary sheaves is outside the present scope.

\subsection{Coordinates with fixed harmonic projection}
Choose an $A_h$-invariant subspace $F_h\subset Q_h\R^{N_h}$ and a mass-orthonormal
basis $V_h$. The two principal choices are the full non-harmonic space, for harmonic
preservation, and $F_h=\Range D^{q-1}$, for cohomology preservation of closed data.
Assume $d_F=\dim F_h>0$. Write
\begin{equation}\label{eq:sector-coordinates}
X=H_hc+V_hz,\quad V_h^\top M_hV_h=I,\quad H_h^\top M_hV_h=0,
\qquad A_F=V_h^\top \mathsf K_hV_h>0.
\end{equation}
The coordinate $c=H_h^\top M_hX$ is fixed along each trajectory. A training law
$\mu_0(\dt c,\dt z)$ can have dependence between $c$ and $z$; score estimation
must retain that dependence through conditioning. If $F_h$ is the exact space,
the data must lie in $H_hc+F_h$. Projecting nonclosed data to that space changes
the modeled data distribution and must be reported as preprocessing.

Let $\mathcal W_t$ be standard Brownian motion in $\R^{r_{\mathrm W}}$, independent of the
initial pair $(c,z_0)$. With deterministic
schedules $\nu(t)>0$ and $b(t)\in\R^{d_F\times r_{\mathrm W}}$, define
\begin{equation}\label{eq:forward-coordinate-sde}
\dt z_t=-\nu(t)A_F z_t\,\dt t+b(t)\,\dt\mathcal W_t,
\qquad a(t)=b(t)b(t)^{\top}.
\end{equation}
Assume $\nu$ and $b$ are continuous on $[0,T]$ and $a(t)$ is positive
definite for every $t>0$. It is then uniformly positive definite on each
$[\varepsilon,T]$, $\varepsilon>0$, and $C_t$ below is invertible for
$t>0$. The transition law is
\begin{align}
E(t,s)&=\exp\!\left[-A_F\int_s^t\nu(r)\,\dt r\right],\qquad
C_t=\int_0^tE(t,s)a(s)E(t,s)^{\top}\,\dt s,\label{eq:forward-covariance}\\
z_t\mid z_0,c&\sim\mathcal N(E(t,0)z_0,C_t).\label{eq:forward-transition}
\end{align}
All gradients and densities below use Lebesgue measure in the $z$ coordinates.
The law on a fixed affine space may be singular in the ambient cochain space,
which causes no problem for this coordinate formulation.

\subsection{Conditional score training and reverse dynamics}
Let $p_t(z\mid c)$ be the conditional forward marginal and
$s_*(z,t,c)=\nabla_z\log p_t(z\mid c)$ its score. The Gaussian transition yields
the denoising-score objective
\begin{equation}\label{eq:score-loss}
\mathcal J(\theta)=\E\!\left[
 w(t)\norm{s_\theta(z_t,t,c)+C_t^{-1}(z_t-E(t,0)z_0)}^2\right],
\end{equation}
where $(c,z_0)\sim\mu_0$, $t$ is sampled independently from a positive
density on $[\varepsilon,T]$, and
$z_t=E(t,0)z_0+C_t^{1/2}\epsilon$ with independent
$\epsilon\sim\mathcal N(0,I)$. Take $\varepsilon>0$ and a positive
bounded weight $w$; training down to $t=0$ requires additional
integrability control.

Conditional expectation of the Gaussian transition score shows that $s_*$ minimizes the population loss over unrestricted square-integrable score fields, when the objective is finite \cite{SongEtAl2021}. 
A chosen parametrized family need not contain this minimizer.

Let $r$ be increasing reverse time and set $t=T-r$. Under the usual
finite-dimensional time-reversal regularity and integrability assumptions
\cite{HaussmannPardoux1986,SongEtAl2021}, the reverse equation on
$0\le r\le T-\varepsilon$ is
\begin{equation}\label{eq:reverse-sde}
\dt\bar z_r=
 [\nu(t)A_F\bar z_r+a(t)s_*(\bar z_r,t,c)]\,\dt r
 +b(t)\,\dt\overline{\mathcal W}_r,\qquad \bar z_0\sim p_T(\cdot\mid c).
\end{equation}
Here $\overline{\mathcal W}_r$ is Brownian motion in the reverse
filtration. The increasing reverse clock determines the signs; no spatial
divergence correction appears because $a$ depends only on time.
Conditional finite Gaussian mixtures with positive component covariances
provide a smooth benchmark. With the exact score and terminal law,
$\bar z_r\sim p_{T-r}(\cdot\mid c)$; stopping at $T-\varepsilon$
yields $p_\varepsilon$, not $p_0$. Recovery at $t=0$ requires the
corresponding endpoint regularity or a limiting argument.

Replacing $s_*$ by $s_\theta$ defines a candidate learned sampler. Its
well-posedness and nonexplosion must be ensured on the sampling interval;
local Lipschitz continuity and at most linear growth of the reverse drift,
uniformly in time for each fixed condition $c$, are sufficient assumptions.
On any interval on which the sampler exists, its reconstruction
\begin{equation}\label{eq:reverse-reconstruct}
\bar X_r=H_hc+V_h\bar z_r
\end{equation}
preserves $P_h\bar X_r=H_hc$ by construction, even if the score is inaccurate.
When $F_h=\Range D^{q-1}$, it also preserves closedness and the discrete
cohomology class. Thus support preservation is an algebraic guarantee of the
chosen coordinates; recovery of the target distribution depends on score accuracy,
terminal initialization, and numerical integration. An unconstrained learned vector
field in ambient coordinates would not have this guarantee.

Under corresponding regularity assumptions, a deterministic probability-flow
alternative uses reverse velocity
\begin{equation}\label{eq:probability-flow}
\frac{\dt\bar z_r}{\dt r}
 =\nu(t)A_F\bar z_r+\tfrac12a(t)s_*(\bar z_r,t,c).
\end{equation}
It has the same positive-time marginals with the exact score and terminal law
\cite{SongEtAl2021}; numerical or learned-score approximations require validation.

\subsection{A concrete noising schedule and the terminal law}
A useful finite-dimensional choice is
\begin{equation}\label{eq:normalized-schedule}
a(t)=2\nu(t)A_F,\qquad C_t=I-E(t,0)^2.
\end{equation}
It leaves $\mathcal N(0,I)$ invariant in the $z$ coordinates and gives an explicit
Gaussian transition. At finite $T$, however, $p_T(\cdot\mid c)$ generally still
depends on the data. If the conditional initial law has finite second moment,
synchronous coupling with a stationary Gaussian gives
\begin{equation}\label{eq:terminal-mismatch}
W_2\bigl(p_T(\cdot\mid c),\mathcal N(0,I)\bigr)
 \le\norm{E(T,0)}_2\,
 W_2\bigl(p_0(\cdot\mid c),\mathcal N(0,I)\bigr),
\end{equation}
where $W_2$ is the quadratic Wasserstein distance and
$\norm{E(T,0)}_2=e^{-\lambda_{\min}(A_F)\int_0^T\nu(r)\,\dt r}$.
Indeed, two coupled trajectories share their noise and differ by $E(T,0)$ times
their initial difference. Infimizing its squared expectation proves the bound.
This quantifies a terminal-prior approximation; it does not make a Gaussian
initialization exact at a finite horizon.

The identity covariance of the invariant law for \eqref{eq:normalized-schedule}
is a finite-dimensional modeling choice. It is not trace class on an infinite-dimensional $L^2$ space and
does not define an $H$-valued stationary continuum field. A continuum-compatible
prior would instead require a trace-class covariance and compatible spectral forcing.

\subsection{Generation and the role of conditioning}
To generate with prescribed harmonic coordinates, fix $c$, initialize the terminal coordinate,
solve the reverse model, and reconstruct using \eqref{eq:reverse-reconstruct}.
To model a distribution over harmonic coordinates, sample $c$ from its empirical or separately
learned law before conditional generation. The constrained diffusion cannot mix
harmonic coordinates on its own. The sets $H_hc+F_h$ are affine spaces of form
coefficients; different values of $c$ do not mean different manifold topologies.

For image applications, a decoder must specify how an evolved form becomes an image.
The preservation results impose no topological property on an arbitrary decoder,
threshold, or segmentation procedure. Data on different manifolds or meshes also
require a correspondence and a common representation before their harmonic
conditions can be compared.

\section{Numerical experiments}
\label{sec:numerical}

Section~\ref{sec:num-annulus} reports deterministic Hodge heat-flow experiments on a planar annulus and a solid torus, including spectral and decay diagnostics.
Section~\ref{sec:num-perturbations} reports a stochastic annulus experiment comparing unprotected, harmonic-projected, and exact noise. 
Further numerical validation and applications are discussed in Section~\ref{sec:discussion}.

All experiments use $q=1$ and a fixed $B\in\{\mathrm{tan},\mathrm{nor}\}$;
the stochastic experiment uses $B=\mathrm{tan}$. We use the matrices and
projector of Section~\ref{sec:discrete}, with $M_h=M_{1,B}$ and
$\mathsf K_h=\mathsf K_{1,B}$, suppressing $B$ where unambiguous.
Diagnostics use mass norms unless a Euclidean norm is indicated.

\subsection{Deterministic diffusion on an annulus and a solid torus}
\label{sec:num-annulus}

\begin{figure}[p]
\centering
\includegraphics[width=\textwidth]{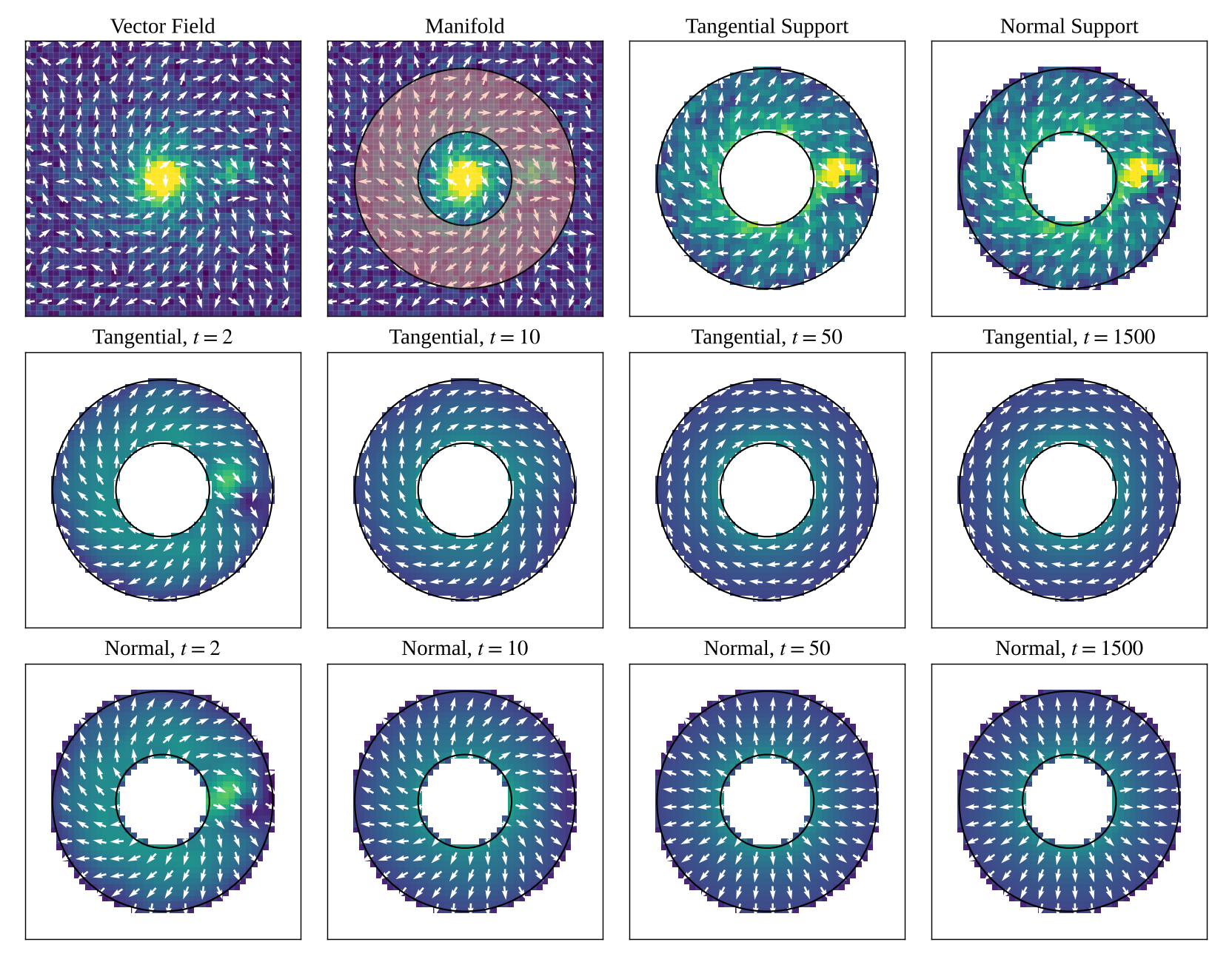}
\caption{\textbf{Deterministic Hodge heat flow of a discrete $1$-form on an
annulus.} Top row, left to right: the initial swirl--radial--bump field with
noise on a $48\times48$ Cartesian grid; the annulus shaded over the field;
and the initial $1$-forms on the tangential and normal supports. Middle and
bottom rows: implicit-Euler evolution with tangential and normal boundary
treatments, respectively, at $t_n=2,10,50,1500$ ($n=4,20,100,3000$). Color
encodes magnitude on a fixed scale within each support row, white arrows
indicate direction, and black curves mark the annular boundary. The
late-time fields are circulating in the tangential case and radial in the
normal case.}
\label{fig:annulus-heat-flow}
\end{figure}

\paragraph{Reported setup and visualization.}
Figure~\ref{fig:annulus-heat-flow} shows a vector field on a $48\times48$
Cartesian grid with unit edge lengths. The initial field combines a swirl,
a radial component, a bump, and noise. An annulus at the center of the grid
is the fixed manifold $M$. The grid field is represented as a discrete
$1$-form on either the tangential or the normal support, and the corresponding
Hodge Laplacian is used for implicit Euler time stepping. The two support-restricted
initial cochains, denoted $X^0_{h,\mathrm{tan}}$ and $X^0_{h,\mathrm{nor}}$, come from the same
underlying grid field but need not be identical vectors or live on identical
sets of cells. The noise enters only through the initial field; no
stochastic forcing is applied during the evolution.

\begin{figure}[p]
\centering
\includegraphics[width=\textwidth]{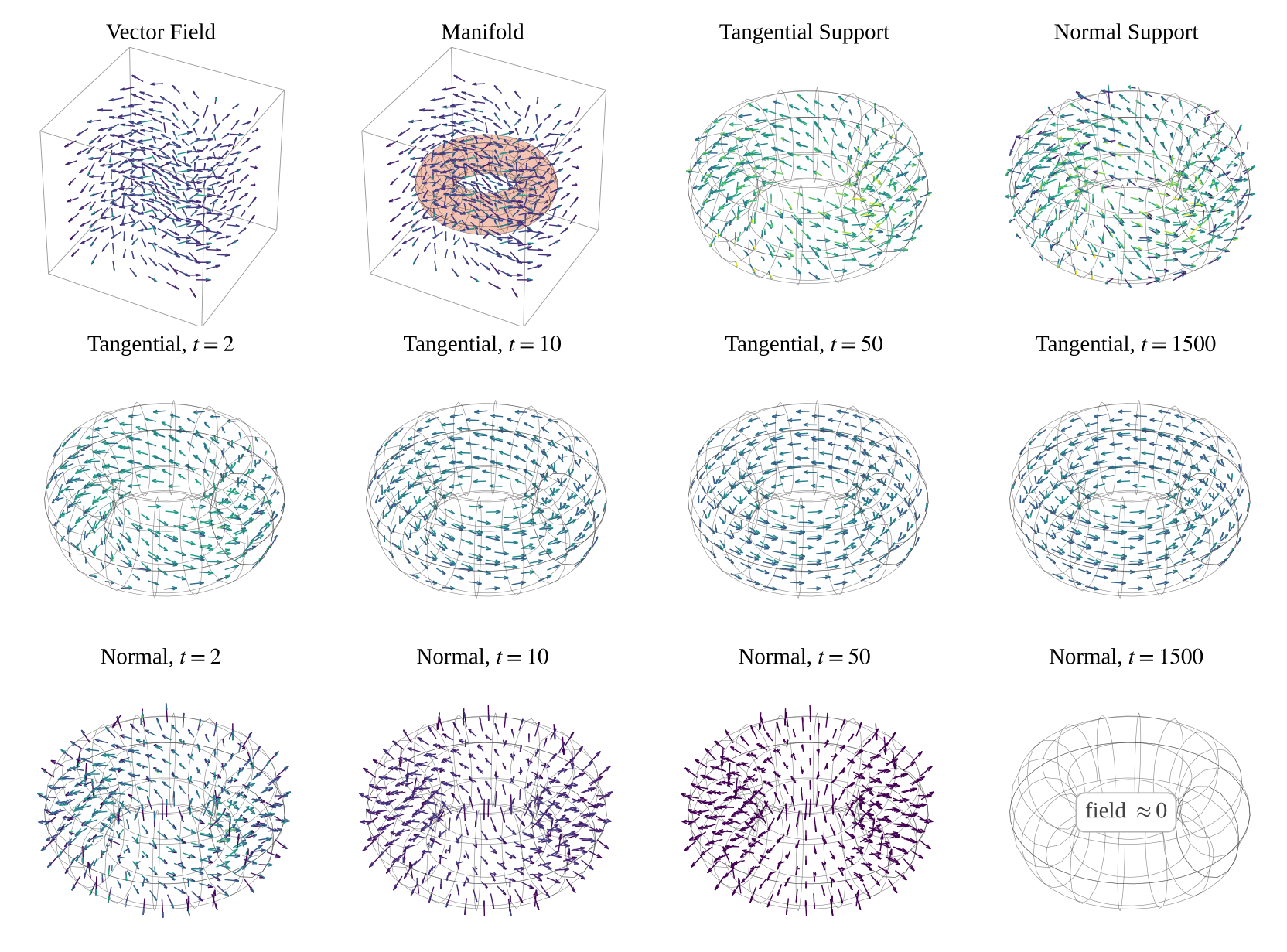}
\caption{\textbf{Deterministic Hodge heat flow of a discrete $1$-form on a
solid torus.} Top row, left to right: the initial swirl--radial--bump field
with noise on a $48\times48\times48$ Cartesian grid, shown at every sixth
cube center; the solid torus shaded over the field; and the initial
$1$-forms on the tangential and normal supports, shown at every third cube
center. Middle and bottom rows: implicit-Euler evolution with tangential and
normal boundary treatments, respectively, at $t_n=2,10,50,1500$
($n=4,20,100,3000$). Arrows have unit length and indicate direction; their
color encodes magnitude on a fixed scale within each support row. The grey
mesh marks $\partial M$. The tangential field relaxes to a circulation
around the central axis. The normal field decays to zero: its maximum
magnitude is $1.8\times10^{-3}$ at $t=50$ and $1.6\times10^{-72}$ at
$t=1500$, where it falls below $10^{-6}$ of the row scale and no arrows
are drawn.}
\label{fig:torus-heat-flow}
\end{figure}

Figure~\ref{fig:torus-heat-flow} repeats the experiment in three dimensions
on a $48\times48\times48$ grid, with a solid torus as $M$. The torus is
chosen so that its cross-section in the plane $z=0$ is exactly the annulus
of the two-dimensional run, and the initial field is the same swirl,
radial, bump, and noise construction extended to three dimensions. 
The same support and operator constructions, time step, and snapshot times are used in both dimensions.

Both runs use $\nu=1$, $\tau=0.5$, and snapshots at
$(n,t_n)=(4,2),(20,10),(100,50),(3000,1500)$. Color represents magnitude
and arrows direction. Each support uses one fixed color scale for its
initial panel and evolution row; the figure captions describe the
remaining display conventions.

\paragraph{Observed behavior.}
On the annulus, the tangential field approaches a circulation around the
hole, while the normal field approaches a radial pattern. On the solid
torus, tangential circulation persists but the entire normal field decays.
Its maximum magnitude falls from $1.36$ initially to $0.63$, $0.17$,
$1.8\times10^{-3}$, and $1.6\times10^{-72}$ at the four displayed times.
Tables~\ref{tab:annulus-diagnostics} and~\ref{tab:torus-diagnostics}
quantify the harmonic preservation and non-harmonic decay behind these
patterns.

\paragraph{Reference fields and continuum interpretation.}
To interpret these patterns using Sections~\ref{sec:hodge} and~\ref{sec:det},
consider the unweighted Euclidean annulus
\[
M=\{(x,y):R_{\mathrm{in}}\le r\le R_{\mathrm{out}}\},
\qquad r=(x^2+y^2)^{1/2},\qquad 0<R_{\mathrm{in}}<R_{\mathrm{out}},
\]
with coordinates centered on the hole at the center of the $N\times N$
grid, $N=48$. The radii are $R_{\mathrm{out}}=0.4N=19.2$ and
$R_{\mathrm{in}}=0.17N=8.16$. 
We identify Euclidean vector fields with their associated $1$-forms using the standard Euclidean metric, and use the same notation for both. 
Define the components of the initial field by
\begin{equation}\label{eq:annulus-reference-harmonics}
 h_{\mathrm{tan}}=\frac{-y\,dx+x\,dy}{r^2}=d\theta,
 \qquad
 h_{\mathrm{nor}}=\frac{x\,dx+y\,dy}{r^2}=d\log r,
 \qquad
 h_b=\nabla_h b,
 \qquad
 h_r = 0.15\, \xi,
\end{equation}
where the total vector field is a swirl, a radial term, a bump, and random
noise, $U= 4h_{\mathrm{tan}}+4h_{\mathrm{nor}}+3h_b+h_r$. The bump is
$b=\exp(-((x-0.28N)^2+y^2)/(0.05N)^2)$, centered at $(x,y)=(13.44,0)$ with
width $2.4$. The discrete gradient $\nabla_h$ uses central differences
inside the grid. The random noise $\xi$ is taken from a zero-mean normal
distribution. Here $d\theta$ denotes the globally defined angular $1$-form, not the
differential of a globally single-valued angle. Direct differentiation
gives $dh_B=\delta h_B=0$ for both fields $h_{\mathrm{tan}}$ and $h_{\mathrm{nor}}$. Their associated vectors are
$r^{-1}e_\theta$ and $r^{-1}e_r$, respectively. Thus circulation need not
mean nonzero curl, and a radial field need not have nonzero divergence on
an annulus whose center is excluded.

The form $h_{\mathrm{tan}}$ satisfies the full tangential/absolute boundary conditions,
while $h_{\mathrm{nor}}$ satisfies the full normal/relative conditions. Both harmonic
spaces have dimension one, so
\begin{equation}\label{eq:annulus-harmonic-spaces}
 \Harm_{\mathrm{tan}}^1=\operatorname{span}\{h_{\mathrm{tan}}\}\simeq H^1(M;\constsheaf),
 \qquad
 \Harm_{\mathrm{nor}}^1=\operatorname{span}\{h_{\mathrm{nor}}\}
       \simeq H^1(M,\partial M;\constsheaf).
\end{equation}
The ordinary class of $h_{\mathrm{tan}}$ is nonzero because its integral along a
positively oriented circle around the hole is $2\pi$. Although $h_{\mathrm{nor}}$ is
exact as an ordinary form, it is not exact in the relative complex: a
primitive $\log r+C$ cannot vanish on both boundary circles. Its integral
along an outward radial path between them is
$\log(R_{\mathrm{out}}/R_{\mathrm{in}})$.
This explains how two one-dimensional harmonic spaces on the same annulus
can have the visibly different directional character in the figure.

For the torus, let $\rho=(x^2+y^2)^{1/2}$ and consider the solid torus
\[
M=\{(x,y,z):(\rho-R_0)^2+z^2\le a^2\},
\qquad R_0=0.285N=13.68,\qquad a=0.115N=5.52,
\]
centered in the $48\times48\times48$ grid, so that $R_0-a=R_{\mathrm{in}}$
and $R_0+a=R_{\mathrm{out}}$. The initial field is
$U=4h_{\mathrm{tan}}+4h_{\mathrm{nor}}+3h_b+h_r$ with
\[
 h_{\mathrm{tan}}=\frac{-y\,dx+x\,dy}{\rho^2}=d\varphi,
 \qquad
 h_{\mathrm{nor}}=\frac{x\,dx+y\,dy}{\rho^2}=d\log\rho,
\]
the bump $b=\exp(-((x-0.28N)^2+y^2+z^2)/(0.05N)^2)$, centered inside the
tube at distance $0.24$ from its core circle, and noise $\xi$ with three
components. Both $h_{\mathrm{tan}}$ and $h_{\mathrm{nor}}$ are again closed and coclosed,
but $h_{\mathrm{nor}}$ does not satisfy the normal boundary condition on the
solid torus. The form
$h_{\mathrm{tan}}$ has vanishing normal component on $\partial M$, because $e_\varphi$ is
tangent to every surface of revolution about the $z$-axis, so it satisfies
the tangential conditions. The solid torus, however, has no normal harmonic
$1$-form:
\begin{equation}\label{eq:torus-harmonic-spaces}
 \Harm_{\mathrm{tan}}^1=\operatorname{span}\{h_{\mathrm{tan}}\},
 \qquad
 \Harm_{\mathrm{nor}}^1=\{0\}.
\end{equation}

For continuum initial data $\omega_{0,B}\in L^2\Omega^1(M)$,
Theorems~\ref{thm:det-invariance} and~\ref{thm:decay} give
\begin{equation}\label{eq:annulus-continuum-limit}
 e^{-\nu t\Delta_B}\omega_{0,B}\longrightarrow P_B\omega_{0,B}
 =\begin{cases}
 \displaystyle\frac{\ip{\omega_{0,B}}{h_B}}{\norm{h_B}^2}h_B,
   &\Harm_B^1=\operatorname{span}\{h_B\},\\[1ex]
 0,&\Harm_B^1=\{0\}.
 \end{cases}
\end{equation}
These are continuum reference calculations: the smooth bump uses $db$,
whereas the implemented field uses the sampled finite-difference gradient
$\nabla_hb$. The reference forms are pointwise orthogonal. Since
$h_{\mathrm{tan}}$ is coclosed with zero normal trace,
$\ip{db}{h_{\mathrm{tan}}}=0$. On the annulus,
$\ip{db}{h_{\mathrm{nor}}}$ is a small boundary contribution, not
identically zero. Thus the corresponding smooth tangential limit is
$4h_{\mathrm{tan}}$ plus the harmonic projection of any $L^2$ noise;
the normal-annulus limit additionally contains the bump contribution.
The discrete projections are measured separately below.

For the solid torus, the tangential reference has
$\norm{h_{\mathrm{tan}}}^2=\int_M\rho^{-2}\,dV
=4\pi^2\bigl(R_0-\sqrt{R_0^2-a^2}\bigr)$.
The normal limit is zero, with decay governed by the smallest excited
positive eigenvalue.

\paragraph{Discrete interpretation and scope of this result.}
For the compatible operators of Section~\ref{sec:discrete}, with
$M_{h,B}$ the diagonal cut-cell Hodge star on edges and
$\mathsf K_{h,B}$ the stiffness-form Hodge Laplacian, each step solves
the implicit-Euler system
\begin{equation}\label{eq:annulus-backward-euler}
 (M_{h,B}+\nu\tau\mathsf K_{h,B})X^{n+1}_{h,B}=M_{h,B}X^{n}_{h,B},
 \quad\text{i.e.}\quad
 X^{n+1}_{h,B}=(I+\nu\tau A_{h,B})^{-1}X^{n}_{h,B},
\end{equation}
with the parameters above. The matrix
$M_{h,B}+\nu\tau\mathsf K_{h,B}$ is factorized once by sparse LU and the
factors are reused at every step, so no iterative solver or stopping
tolerance is involved; the relative residual of a solve on the tangential
annulus is $6\times10^{-16}$. By Section~\ref{sec:discrete}, the update
fixes $P_{h,B}X^0_{h,B}$ in exact arithmetic and multiplies each positive
mode by $(1+\nu\tau\lambda)^{-1}$. An explicit step at the same $\tau$
would violate \eqref{eq:em-cfl}: $\nu\tau\lambda_{h,B,\max}$ is approximately
$16.4$ (tangential) and $25.8$ (normal) on the annulus, and $1192$
(tangential) and $282$ (normal) on the torus. The large torus values come
from small cut cells. The smallest tangential Hodge-star entry, a cut dual
face area, is $6.4\times10^{-6}$, and the largest normal entry, the
reciprocal of a cut edge length, is $207$.

The initial cochains are not closed. With $D^1$ the support-restricted
discrete exterior derivative on $1$-cochains, the unweighted Euclidean norm
$\norm{D^1X^0_{h,B}}_2$ on the
annulus is $6.4$ (tangential) and $7.5$ (normal), compared with
$\norm{X^0_{h,B}}_2\approx15$. Most of this comes from the noise. Without
it, the residual falls to $0.19$ on the tangential support, which is the
discretization error from sampling the closed continuum terms. On the
normal support it remains $4.0$, almost entirely on cut boundary faces,
where the cut-edge scaling of the restriction does not preserve
closedness. The torus shows the same pattern. There
$\norm{D^1X^0_{h,B}}_2$ is $34.3$ (tangential) and $39.0$
(normal), compared with $\norm{X^0_{h,B}}_2$ of $46.5$ and $42.0$. Without
the noise it falls to $0.57$ on the tangential support, almost all on whole
faces, and to $21.5$ on the normal support, of which $21.46$ lies on cut
faces. Preservation of a cohomology class of the initial field is
therefore not tested. What is tested is preservation of the harmonic
projection $P_{h,B}X^0_{h,B}$.

\paragraph{Quantitative diagnostics.}
The diagnostics
\begin{align}
 \epsilon_{H,B}^{n}
   &=\norm{P_{h,B}(X^n_{h,B}-X^0_{h,B})}_{M_{h,B}},
       \label{eq:annulus-harmonic-diagnostic}\\
 E_{\perp,B}^{n}
   &=\tfrac12\norm{(I-P_{h,B})X^n_{h,B}}_{M_{h,B}}^2,
 \qquad
 \mathcal R_B^{n}=\norm{A_{h,B}X^n_{h,B}}_{M_{h,B}},
       \label{eq:annulus-decay-diagnostic}
\end{align}
test invariant preservation, non-harmonic relaxation, and approach to a
discrete harmonic field separately. Here $P_{h,B}$ is the
$M_{h,B}$-orthogonal projector onto the discrete harmonic space, built from
the smallest eigenpairs of $\mathsf K_{h,B}h=\lambda M_{h,B}h$ computed by
shift-invert Lanczos (ARPACK), with shift $-10^{-6}$ and the default
machine-precision convergence tolerance. An eigenvalue is identified as
zero when $\lambda\le10^{-10}$, and the corresponding eigenvectors are
$M_{h,B}$-orthonormalized. In each one-dimensional case, $H_{h,B}$ denotes
the resulting $M_{h,B}$-normalized basis column. When the harmonic space is trivial, $P_{h,B}=0$,
$\epsilon_{H,B}^n$ is identically zero, and $E_{\perp,B}^n$ is the full
energy.

On the annulus the harmonic basis was checked against a dense generalized
eigendecomposition. On each support exactly one eigenvalue lies below
$10^{-10}$ (in magnitude below $10^{-15}$), so the discrete harmonic
dimension is one in both cases, matching \eqref{eq:annulus-harmonic-spaces}.
The two computations give the same kernel vector up to sign to machine
precision, and the residual $\norm{A_{h,B}H_{h,B}}_{M_{h,B}}$ is below
$10^{-15}$. 
The first positive eigenvalue is $\lambda_{h,\mathrm{tan},1}=5.587\times10^{-3}$ and $\lambda_{h,\mathrm{nor},1}=5.580\times10^{-3}$, each with an observed twofold multiplicity, consistent with the symmetries of the discretized annulus. 
The zero cutoff is therefore not delicate: the computed zero eigenvalues
are separated from the first positive eigenvalue by more than twelve
orders of magnitude, here and on the torus.
The largest eigenvalues are $32.86$ and $51.60$. The supports carry
$1992$ (tangential) and $2012$ (normal) edge unknowns.

The torus supports carry $28{,}924$ (tangential) and $26{,}832$ (normal)
edge unknowns, too many for a dense eigendecomposition, so the basis there
is checked through its residual and the spectral gap. On the tangential
support exactly one eigenvalue lies below $10^{-10}$ (it is
$2.7\times10^{-16}$), with residual $\norm{A_{h,\mathrm{tan}}H_{h,\mathrm{tan}}}_{M_{h,\mathrm{tan}}}=2.8\times10^{-15}$.
The next eigenvalue is $\lambda_{h,\mathrm{tan},1}=5.540\times10^{-3}$, with
multiplicity two. On the normal support no eigenvalue lies below
$10^{-10}$. The smallest are $0.10956$, $0.11271$, and $0.11525$ (the last
with multiplicity two). The discrete harmonic dimensions are therefore one
and zero, matching \eqref{eq:torus-harmonic-spaces}. The largest
eigenvalues are $2385$ and $564$.

\begin{table}[ht]
\centering
\small
\begin{tabular}{rr|ccc|ccc}
\hline
 & & \multicolumn{3}{c|}{tangential} & \multicolumn{3}{c}{normal}\\
$n$ & $t_n$ & $\epsilon_{H,\mathrm{tan}}^n$ & $E_{\perp,\mathrm{tan}}^n$ & $\mathcal R_{\mathrm{tan}}^n$
           & $\epsilon_{H,\mathrm{nor}}^n$ & $E_{\perp,\mathrm{nor}}^n$ & $\mathcal R_{\mathrm{nor}}^n$\\
\hline
0    & 0    & $0$                  & $6.60\times10^{1}$  & $3.51\times10^{1}$
            & $0$                  & $6.56\times10^{1}$  & $2.94\times10^{1}$\\
4    & 2    & $1.0\times10^{-15}$  & $2.81\times10^{1}$  & $1.40$
            & $1.2\times10^{-16}$  & $2.90\times10^{1}$  & $1.32$\\
20   & 10   & $6.4\times10^{-15}$  & $6.72$              & $3.14\times10^{-1}$
            & $3.6\times10^{-15}$  & $7.02$              & $3.25\times10^{-1}$\\
100  & 50   & $2.7\times10^{-14}$  & $3.16\times10^{-2}$ & $1.14\times10^{-2}$
            & $2.6\times10^{-14}$  & $1.74\times10^{-2}$ & $1.15\times10^{-2}$\\
3000 & 1500 & $8.2\times10^{-13}$  & $1.88\times10^{-9}$ & $3.43\times10^{-7}$
            & $9.1\times10^{-13}$  & $7.00\times10^{-10}$& $2.09\times10^{-7}$\\
\hline
\end{tabular}
\caption{Diagnostics \eqref{eq:annulus-harmonic-diagnostic}--\eqref{eq:annulus-decay-diagnostic}
for the annulus run in Figure~\ref{fig:annulus-heat-flow}.}
\label{tab:annulus-diagnostics}
\end{table}

\begin{table}[ht]
\centering
\small
\begin{tabular}{rr|ccc|cc}
\hline
 & & \multicolumn{3}{c|}{tangential} & \multicolumn{2}{c}{normal}\\
$n$ & $t_n$ & $\epsilon_{H,\mathrm{tan}}^n$ & $E_{\perp,\mathrm{tan}}^n$ & $\mathcal R_{\mathrm{tan}}^n$
           & $E_{\perp,\mathrm{nor}}^n$ & $\mathcal R_{\mathrm{nor}}^n$\\
\hline
0    & 0    & $0$                  & $5.46\times10^{2}$  & $1.30\times10^{2}$
            & $9.06\times10^{2}$   & $8.78\times10^{2}$\\
4    & 2    & $2.8\times10^{-14}$  & $1.95\times10^{2}$  & $3.65$
            & $3.24\times10^{2}$   & $5.03$\\
20   & 10   & $1.3\times10^{-13}$  & $3.07\times10^{1}$  & $9.10\times10^{-1}$
            & $4.02\times10^{1}$   & $1.15$\\
100  & 50   & $5.8\times10^{-13}$  & $7.36\times10^{-2}$ & $1.14\times10^{-2}$
            & $5.27\times10^{-3}$  & $1.16\times10^{-2}$\\
3000 & 1500 & $1.9\times10^{-11}$  & $6.07\times10^{-9}$ & $6.11\times10^{-7}$
            & $4.32\times10^{-141}$& $1.05\times10^{-71}$\\
\hline
\end{tabular}
\caption{Diagnostics \eqref{eq:annulus-harmonic-diagnostic}--\eqref{eq:annulus-decay-diagnostic}
for the torus run in Figure~\ref{fig:torus-heat-flow}. The normal harmonic
space is trivial, so $\epsilon_{H,\mathrm{nor}}^n\equiv0$ is omitted and
$E_{\perp,\mathrm{nor}}^n$ is the full energy.}
\label{tab:torus-diagnostics}
\end{table}

Tables~\ref{tab:annulus-diagnostics} and~\ref{tab:torus-diagnostics} show
the three behaviors separately. First, the harmonic drift
$\epsilon_{H,B}^n$ stays at round-off level, below $10^{-12}$ on the
annulus and below $2\times10^{-11}$ on the torus after $3000$ steps,
consistent with exact preservation of $P_{h,B}X^0_{h,B}$ by
\eqref{eq:annulus-backward-euler}. The harmonic coefficients have
magnitudes $9.4099$ (annulus, tangential), $9.4497$ (annulus, normal), and
$26.7479$ (torus, tangential); their signs depend only on the sign of the
computed basis vector. Second, $E_{\perp,B}^n$ decreases monotonically. On
the annulus it falls by factors of $3.5\times10^{10}$ and $9.4\times10^{10}$,
and its late-time per-step ratio, $0.99443592$ (tangential) and
$0.99444349$ (normal), agrees with the predicted slowest-mode factor
$(1+\nu\tau\lambda_{h,B,1})^{-2}$ to about $10^{-10}$. On the torus the
tangential energy falls by a factor of $9.0\times10^{10}$, with late-time
ratio $0.99448270$, again matching $(1+\nu\tau\lambda_{h,\mathrm{tan},1})^{-2}$ to
about $10^{-9}$. The normal energy on the torus falls by $143$ orders of
magnitude. 
The Laplacian residuals fall to order $10^{-7}$ in the three runs with
a harmonic field and to order $10^{-71}$ in the normal-torus run. The
distance to the discrete harmonic space is $\sqrt{2E_{\perp,B}^n}$,
not the Laplacian residual; the two are related by
$\norm{(I-P_{h,B})X^n_{h,B}}_{M_{h,B}}\le\mathcal R_B^n/\lambda_{h,B,1}$.

The normal decay rate on the torus needs one further remark. Over the
window $n=400\to401$ ($t=200\to200.5$), the per-step energy ratio
$E_{\perp,\mathrm{nor}}^{n+1}/E_{\perp,\mathrm{nor}}^{n}$ is $0.89614538$. This is close to the
reported factor $0.89614561$ associated with the second eigenvalue
$\lambda\approx0.11271$, rather than $0.89882750$ associated with the
smallest eigenvalue $\lambda\approx0.10956$. Eigenvalues are rounded
here; high-precision ratio comparisons require their unrounded values. Projecting $X^0_{h,\mathrm{nor}}$ onto the
$M_{h,\mathrm{nor}}$-orthonormal slowest eigenvectors explains this. The coefficient
on the second mode is $-24.672$, almost entirely from the radial term
$4h_{\mathrm{nor}}$. The coefficient on the first mode is $-0.055$ and comes only from
the noise; the radial term contributes nothing above round-off there.
With $c_k$ these coefficients and $g_k=(1+\nu\tau\lambda_k)^{-1}$, the
energy is $E_{\perp,\mathrm{nor}}^n=\tfrac12\sum_k c_k^2g_k^{2n}$. At $n=400$ the two
slowest modes alone predict a ratio of $0.89614565$. Adding the next four
modes, with eigenvalues $0.11525$ and $0.11763$ (each double, so only the
combined weight of each pair is basis independent), gives $0.89614538$,
equal to the measured value to all displayed digits. The ratio is not yet
constant at $n=400$. It drifts upward as the first mode gains weight, to
$0.89614583$ at $n=1000$, $0.89615088$ at $n=2000$, and $0.89624653$ at
$n=3000$, and the six-mode prediction matches the measurement to
$10^{-8}$ in each window. The two slowest modes carry equal energy at
$n\approx4085$; only beyond the displayed range would the ratio approach
the first-mode factor.

The harmonic coefficients can also be compared with
\eqref{eq:annulus-continuum-limit}. Splitting $X^0_{h,B}$ by the terms of
$U$ shows the following, with signs taken relative to the total
coefficient. On the annulus, in the tangential case $4h_{\mathrm{tan}}$
contributes $9.287$ of the coefficient, $4h_{\mathrm{nor}}$ and $3h_b$ contribute zero to
displayed precision, and the noise contributes $0.123$. In the normal case
$4h_{\mathrm{nor}}$ contributes $9.280$, $4h_{\mathrm{tan}}$ contributes zero, $3h_b$ contributes
$-0.005$, and the noise contributes $0.175$. The continuum value
$4\norm{h_B}=4\,(2\pi\log(R_{\mathrm{out}}/R_{\mathrm{in}}))^{1/2}\approx9.275$,
identical for both $B$, is within $0.15\%$ of the discrete swirl and radial
contributions. Relative to the normalized continuum coefficients, the measured totals
have ratios approximately $1.015$ and $1.019$, mainly reflecting the
harmonic contribution of the initialization noise. On the torus, $4h_{\mathrm{tan}}$ contributes $26.936$ of the tangential
coefficient, $4h_{\mathrm{nor}}$ and $3h_b$ contribute zero to round-off, and the noise
contributes $0.188$ with the opposite sign. The continuum value
$4\norm{h_{\mathrm{tan}}}\approx27.105$ is $0.6\%$ above the discrete swirl
contribution on this mesh, whose tube radius is only about $5.5$ cells.
The total coefficient is approximately $0.987$ times its continuum
reference. These coefficient comparisons do not measure continuum field
error or establish a mesh-convergence rate.

For reproducibility, the remaining parameters are as follows. The annulus
is the zero sublevel set of
$\max(r-R_{\mathrm{out}},\,R_{\mathrm{in}}-r)$ sampled at grid vertices,
with center $(23.5,23.5)$ in vertex coordinates. The torus is the zero
sublevel set of $((\rho-R_0)^2+z^2)^{1/2}-a$ sampled at grid vertices, with
center $(23.5,23.5,23.5)$. In both cases $r^2$ or $\rho^2$ is floored at $1$
in the initial field, which affects only points inside the hole. The noise
$\xi$ is i.i.d.\ standard normal in each component at each vertex, with
amplitude $0.15$ and NumPy seed $0$, and is added once to the initial
field. The vector-to-cochain map assigns to each edge the average of the
relevant component at its two endpoints. On the normal support this value
is multiplied by the cut primal edge length. Displayed vectors average the
cochain values on the parallel edges of each face (annulus) or cube
(torus), with the normal-support scaling undone first. The tangential
support consists of the edges whose dual cell meets $M$: the dual edge in
two dimensions and the dual face, the unit square through four cube
centers, in three dimensions. $M_{h,\mathrm{tan}}$ is the cut dual edge length or the
cut dual face area, used as computed with no lower bound. The normal
support consists of the edges with an endpoint vertex inside $M$, with
$M_{h,\mathrm{nor}}$ the reciprocal of the cut primal edge length. The cut
length is bounded below by $10^{-5}$ before the reciprocal is taken, so
the entries of $M_{h,\mathrm{nor}}$ are at most $10^{5}$. In addition,
sampled level-set values within $10^{-5}$ of zero are moved to
$\pm10^{-5}$, with an exact zero counted as outside $M$. Neither bound was
active in these runs: the shortest cut edge lengths are $0.073$ (annulus)
and $0.0048$ (torus), and the smallest sampled level-set magnitudes are
$0.014$ and $5.9\times10^{-5}$.
Convergence under grid refinement was not studied for either domain.
\FloatBarrier

\subsection{Stochastic diffusion on an annulus}
\label{sec:num-perturbations}

We report a stochastic experiment on the tangential (absolute) annulus complex. 
Starting from the same closed initial cochain, we compare unprotected, harmonic-projected, and exact noise, using identical Brownian increments across the three runs.
We measure harmonic drift and the closedness defect to distinguish harmonic-coordinate preservation from cohomology preservation.

\begin{figure}[t]
\centering
\includegraphics[width=\textwidth]{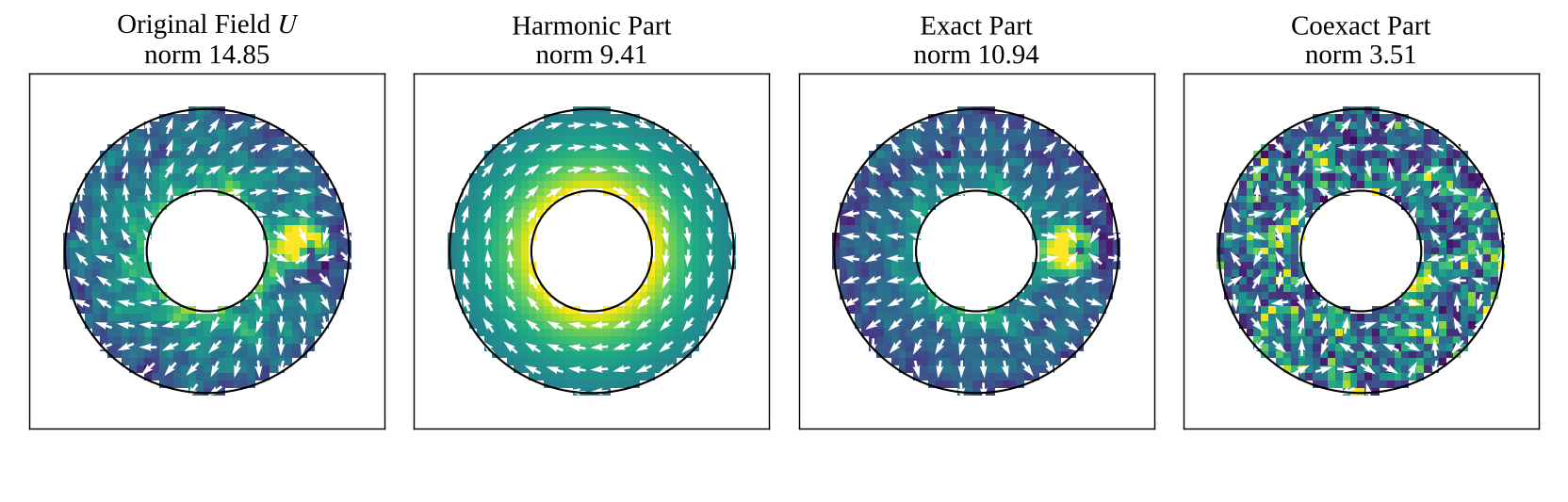}
\caption{\textbf{Hodge decomposition used for the stochastic example.} Left
to right: the tangential cochain $U_h$ of the initial field of
Section~\ref{sec:num-annulus}, and its harmonic, exact, and coexact parts
\eqref{eq:stoch-hodge-split}, with their $M_h$-norms. Color encodes
magnitude on a separate scale in each panel, white arrows indicate
direction, and black curves mark $\partial M$. The normalized parts are the
noise directions $\eta_h$, $e_h$, and $v_h$.}
\label{fig:stoch-setup}
\end{figure}

\paragraph{Setup.}
We reuse the $1992$-edge tangential annulus complex from
Section~\ref{sec:num-annulus}. Noise directions are obtained by decomposing
its initial cochain $U_h=X^0_{h,\mathrm{tan}}$ as
\begin{equation}\label{eq:stoch-hodge-split}
 U_h=U_{\mathrm H}+U_{\mathrm E}+U_{\mathrm C},
 \qquad
 U_{\mathrm H}=P_hU_h,
 \qquad
 U_{\mathrm E}=D^0\phi_h,
 \qquad
 U_{\mathrm C}=U_h-U_{\mathrm H}-U_{\mathrm E},
\end{equation}
where $\phi_h$ solves the discrete Poisson problem
$(D^0)^\top M_hD^0\phi_h=(D^0)^\top M_hU_h$, so that $U_{\mathrm E}$ is the
$M_h$-orthogonal projection of $U_h$ onto the exact cochains and the
remainder $U_{\mathrm C}$ is coexact. The scalar potential is determined
only modulo constants; its gradient $D^0\phi_h$ is unique. Numerically,
the one-dimensional constant nullspace is removed by a small mass shift:
we solve with $(D^0)^\top M_hD^0+10^{-12}M_0$, where $M_0$ is the vertex
Hodge star, by sparse LU. The right-hand side is compatible, with entries
summing to $4\times10^{-15}$, and the resulting $U_{\mathrm E}$ agrees with
the one obtained from a dense pseudoinverse solve to a relative
$10^{-11}$.
The parts have norms
\[
 \norm{U_h}_{M_h}=14.85,\qquad
 \norm{U_{\mathrm H}}_{M_h}=9.41,\qquad
 \norm{U_{\mathrm E}}_{M_h}=10.94,\qquad
 \norm{U_{\mathrm C}}_{M_h}=3.51,
\]
and are shown in Figure~\ref{fig:stoch-setup}. The harmonic component
has coefficient $9.4099$. The exact part contains the projected radial and
bump contributions and the exact part of the initialization noise. The
coexact part is dominated by the noise, with sampling residuals as reported
in Section~\ref{sec:num-annulus}. Normalizing gives
\[
 \eta_h=\frac{U_{\mathrm H}}{\norm{U_{\mathrm H}}_{M_h}},\qquad
 e_h=\frac{U_{\mathrm E}}{\norm{U_{\mathrm E}}_{M_h}},\qquad
 v_h=\frac{U_{\mathrm C}}{\norm{U_{\mathrm C}}_{M_h}},\qquad
 X^0=\eta_h+e_h .
\]
The assumptions used below were checked on these arrays. The Gram matrix of
$[\eta_h,e_h,v_h]$ in the $M_h$ inner product differs from the identity by
less than $3\times10^{-11}$. The restricted operators satisfy $D^1D^0=0$
exactly. The curls $\norm{D^1\eta_h}_{M_2}$ and $\norm{D^1e_h}_{M_2}$ are
below $10^{-15}$, and $\norm{D^1v_h}_{M_2}=1.817$. The initial field $X^0$ is
therefore closed to numerical precision, with harmonic coefficient
$c^0=\eta_h^\top M_hX^0=1$ to the reported precision. The predictions
below use the exact identities of the compatible complex.

The normalized coexact direction has Rayleigh quotient $v_h^\top\mathsf K_hv_h=3.30$, corresponding to an initial squared-norm decay timescale of approximately $0.15$.
Since $v_h$ is not an eigenmode, its evolution involves multiple decay rates. 
We therefore compare the measured statistics with the discrete prediction \eqref{eq:stoch-closedness-prediction}.

\paragraph{Noise constructions and update.}
Three noise matrices are compared, corresponding to unprotected,
harmonic-projected, and exact noise:
\begin{equation}\label{eq:stoch-noise-matrices}
 B_h^{(1)}=\sigma[\eta_h,e_h,v_h],\qquad
 B_h^{(2)}=\sigma[0,e_h,v_h],\qquad
 B_h^{(3)}=\sigma[0,e_h,0].
\end{equation}
Here $B_h^{(2)}=(I-P_h)B_h^{(1)}$, and $B_h^{(3)}$ additionally removes the
coexact column. Each run uses the implicit update
\begin{equation}\label{eq:stoch-update}
 (M_h+\nu\tau\mathsf K_h)X^{n+1}
 =M_hX^n+\sqrt{\tau}\,M_hB_h\xi_n,
 \qquad \xi_n\sim\mathcal N(0,I_3),
\end{equation}
that is, $X^{n+1}=R_h(X^n+\sqrt{\tau}\,B_h\xi_n)$ with
$R_h=(I+\nu\tau A_h)^{-1}$. As in Section~\ref{sec:num-annulus}, the
system matrix is factorized once by sparse LU and the factors are reused
for every step and trajectory. The increments are independent across steps and
trajectories. The parameters are $\nu=1$, $\tau=0.5$, $T=50$ ($100$ steps),
$\sigma=0.1$, and $200$ independent trajectories generated with NumPy seed
$0$. The three runs use the same increments $\xi_n$, so run~2 is run~1
without its harmonic kicks and run~3 is run~2 without its coexact kicks.
The injected powers $\operatorname{tr}(B_h^\top M_hB_h)$ are $3\sigma^2$,
$2\sigma^2$, and $\sigma^2$. Columns are removed without rescaling the
remaining ones, so this example deliberately does not match the total
injection; it isolates the effect of each removed component.

\paragraph{Predictions for the discrete scheme.}
Unrolling \eqref{eq:stoch-update} gives
\begin{equation}\label{eq:stoch-unrolled}
 X^n=R_h^nX^0+\sqrt{\tau}\sum_{m=0}^{n-1}R_h^{\,n-m}B_h\xi_m .
\end{equation}
Two second moments follow without sampling. First, because
$\eta_h^\top M_hR_h=\eta_h^\top M_h$, the harmonic coefficient
$c^n=\eta_h^\top M_hX^n$ satisfies
$c^{n+1}=c^n+\sqrt{\tau}\,\eta_h^\top M_hB_h\xi_n$. In run~1 this is
$c^n-c^0=\sigma\sqrt{\tau}\sum_{m<n}\xi_{m,1}$, a Gaussian random walk, so
\begin{equation}\label{eq:stoch-drift-prediction}
 c^n-c^0\sim\mathcal N(0,\sigma^2t_n),
 \qquad
 \mathbb{E}\,\norm{P_h(X^n-X^0)}_{M_h}^2=\sigma^2t_n.
\end{equation}
This agrees with \eqref{eq:exact-kernel-noise-leakage}, taking
$H_h=\eta_h$ so that $H_h^\top M_hB_h^{(1)}=\sigma[1,0,0]$. It holds exactly for the scheme, not
only as $\tau\to0$. In runs~2 and~3 the coefficient is constant. Second,
$D^1R_h^nX^0=0$ and $D^1R_h^j\eta_h=D^1R_h^je_h=0$, because the update maps
closed cochains to closed cochains. Only the coexact column contributes to
the curl, and independent increments add their variances:
\begin{equation}\label{eq:stoch-closedness-prediction}
 \mathbb{E}\,\norm{D^1X^n}_{M_2}^2
 =\sigma^2\tau\sum_{j=1}^{n}\norm{D^1R_h^{\,j}v_h}_{M_2}^2
\end{equation}
in runs~1 and~2, and zero in run~3. The terms decay geometrically in $j$,
so the right-hand side rises to a plateau, the balance between injected and
dissipated curl, in contrast to the unbounded growth in
\eqref{eq:stoch-drift-prediction}. For the present $v_h$ the sum is
evaluated numerically; it reaches $89\%$ of its limit after two steps and
its limit is $2.915\times10^{-3}$.
This is the stationary mean squared curl of the implicit scheme at $\tau=0.5$. 
For the corresponding continuous-time semidiscrete SDE, the unit-norm coexact forcing gives stationary mean squared curl $\sigma^2/(2\nu)=0.005$.
The smaller discrete value reflects the covariance bias of the implicit time step.

\begin{figure}[p]
\centering
\includegraphics[width=\textwidth]{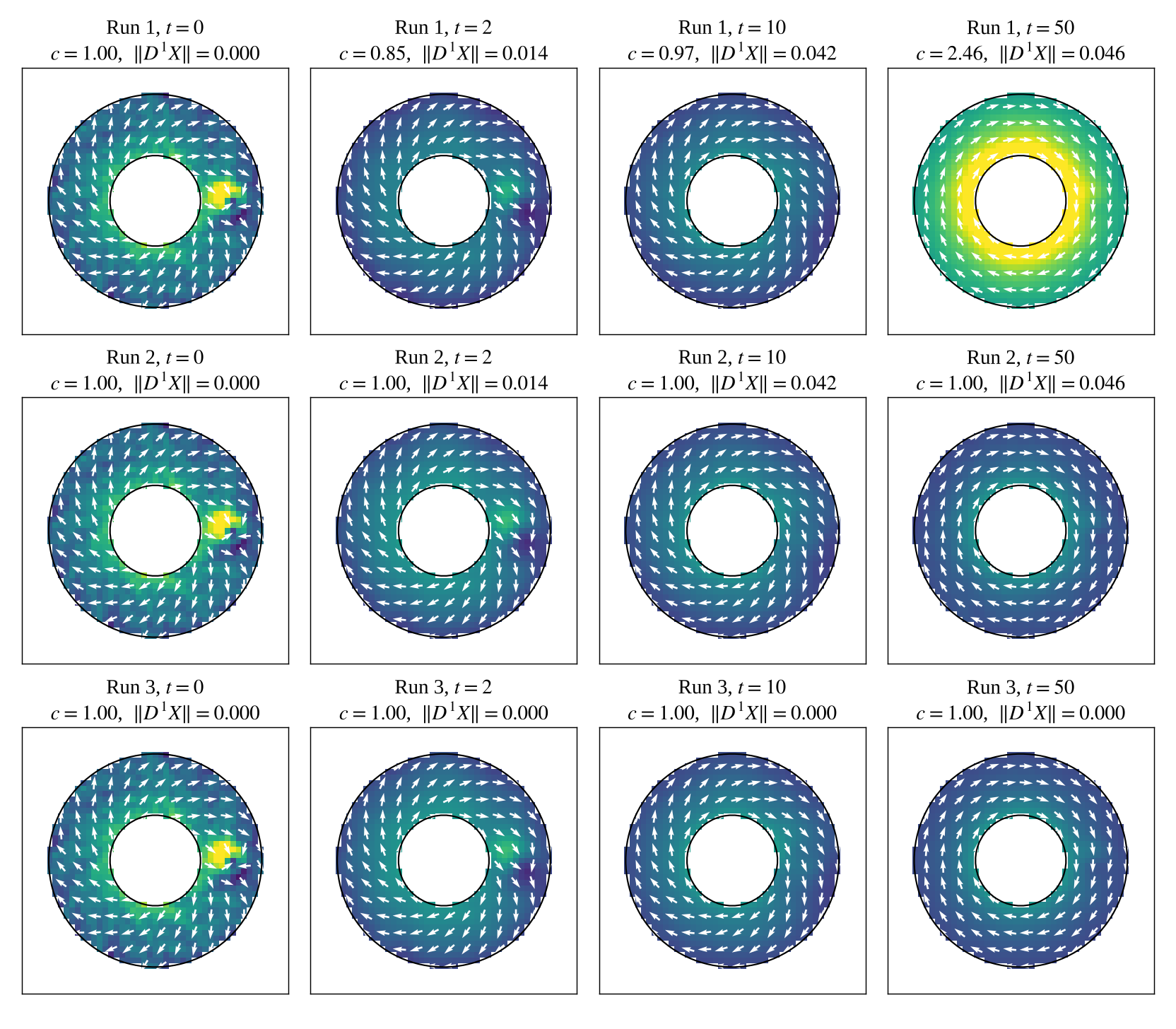}
\caption{\textbf{One sample path of each stochastic run.} Rows: run~1
(unprotected noise $\sigma[\eta_h,e_h,v_h]$), run~2 (harmonic-projected
noise $\sigma[0,e_h,v_h]$), and run~3 (exact noise $\sigma[0,e_h,0]$), all
from $X^0=\eta_h+e_h$ with the same Brownian increments. Columns: $t=0$,
$2$, $10$, $50$. Each title gives that path's harmonic coefficient
$c=\eta_h^\top M_hX^n$ and curl $\norm{D^1X^n}_{M_2}$. Color encodes
magnitude on one scale for all panels, white arrows indicate direction, and
black curves mark $\partial M$.}
\label{fig:stoch-snapshots}
\end{figure}

\begin{figure}[t]
\centering
\includegraphics[width=\textwidth]{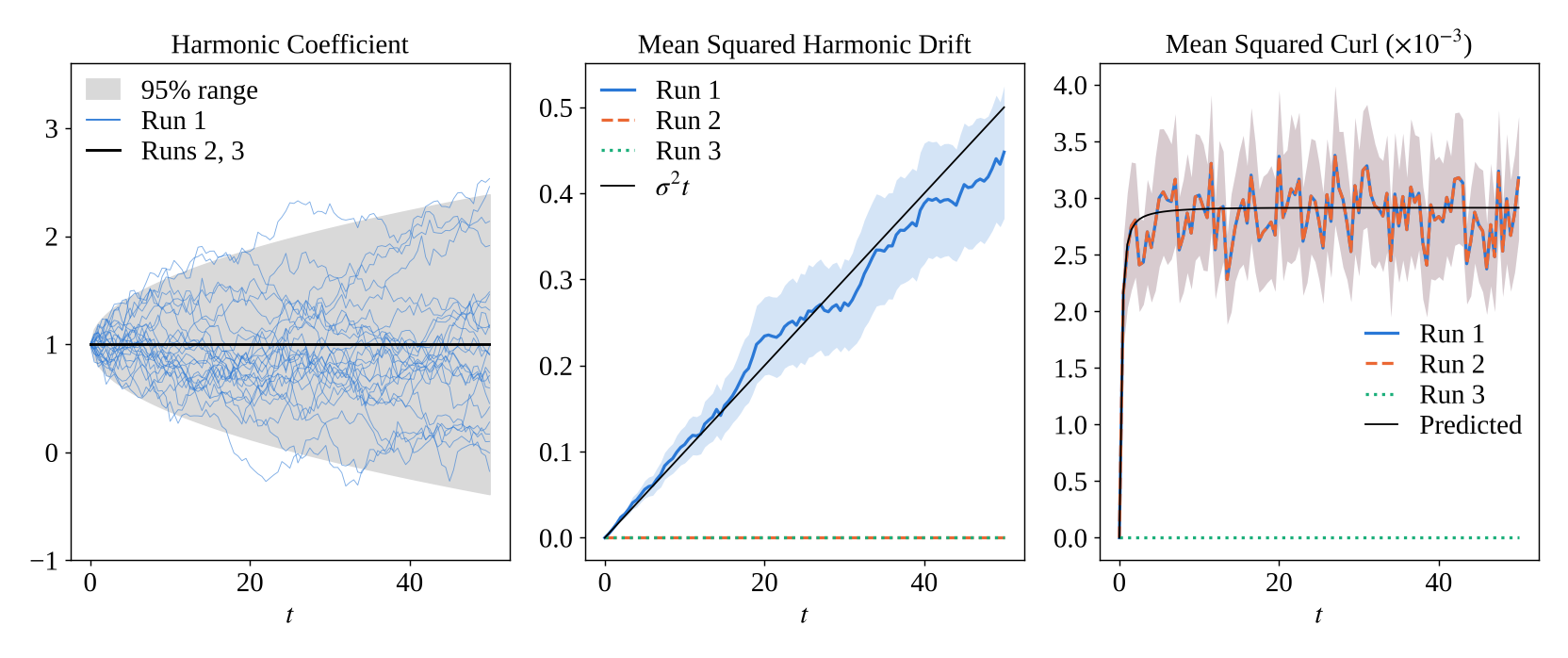}
\caption{\textbf{Harmonic drift and closedness in the stochastic example.}
Left: the harmonic coefficient $c^n$ along $25$ trajectories of run~1,
with the pointwise $95\%$ Gaussian interval
$1\pm1.96\,\sigma\sqrt{t}$; this is not a simultaneous pathwise band.
Runs~2 and~3 remain at $c=1$ to numerical precision. Center:
sample mean over $200$ trajectories of
$\norm{P_h(X^n-X^0)}_{M_h}^2$ with the prediction $\sigma^2t$. Right: sample
mean of $\norm{D^1X^n}_{M_2}^2$, in units of $10^{-3}$, with the prediction
\eqref{eq:stoch-closedness-prediction}; the curves of runs~1 and~2
coincide. Shaded bands are pointwise $95\%$ confidence intervals for
the ensemble means.}
\label{fig:stoch-diagnostics}
\end{figure}

\paragraph{Results.}
Figure~\ref{fig:stoch-snapshots} shows one trajectory of each run and
Figure~\ref{fig:stoch-diagnostics} the ensemble diagnostics;
Table~\ref{tab:stoch-diagnostics} lists the sample means with $95\%$
confidence half-widths. The confidence intervals use the normal
approximation, $\bar x\pm1.96\,s/\sqrt{200}$, where $\bar x$ and $s$ are
the sample mean and sample standard deviation over the $200$ trajectories,
computed separately at each step; no bootstrap was used.

\begin{table}[ht]
\centering
\small
\begin{tabular}{rr|cc|cc}
\hline
 & & \multicolumn{2}{c|}{$\norm{P_h(X^n-X^0)}_{M_h}^2$, run 1}
   & \multicolumn{2}{c}{$10^{3}\,\norm{D^1X^n}_{M_2}^2$, runs 1 and 2}\\
$n$ & $t_n$ & sample mean & $\sigma^2t_n$ & sample mean & predicted\\
\hline
4   & 2  & $0.024\pm0.004$ & $0.020$ & $2.81\pm0.50$ & $2.78$\\
20  & 10 & $0.108\pm0.021$ & $0.100$ & $3.03\pm0.53$ & $2.90$\\
50  & 25 & $0.253\pm0.052$ & $0.250$ & $2.80\pm0.52$ & $2.91$\\
100 & 50 & $0.448\pm0.077$ & $0.500$ & $3.18\pm0.54$ & $2.91$\\
\hline
\end{tabular}
\caption{Sample means over $200$ trajectories, with pointwise $95\%$
confidence half-widths, against the predictions
\eqref{eq:stoch-drift-prediction} and
\eqref{eq:stoch-closedness-prediction}. In runs~2 and~3 the squared
harmonic drift is below $3\times10^{-29}$ on every trajectory and step, and
in run~3 the squared curl is below $4\times10^{-31}$.}
\label{tab:stoch-diagnostics}
\end{table}

The ensemble statistics are consistent with the discrete predictions.
At $T=50$, run~1 has mean squared harmonic drift $0.448\pm0.077$
against $\sigma^2T=0.5$, and coefficient standard deviation $0.666$
against $\sigma\sqrt T=0.707$. The coefficient lies in the pointwise
Gaussian interval for $96.5\%$ of the sampled paths. Runs~2 and~3
preserve it to the numerical precision reported in
Table~\ref{tab:stoch-diagnostics}.

The mean squared curl in runs~1 and~2 approaches the discrete plateau.
Their curls agree path by path to $4\times10^{-17}$ because the paired
runs differ only by harmonic increments. Run~3 preserves closedness.
Together these diagnostics distinguish harmonic-coordinate preservation
(run~2) from preservation of the initial closed cohomology class (run~3).
The displayed run~1 trajectory ends at $c=2.46$, about $2.1$
standard deviations above its initial value. Runs~2 and~3 appear similar
in the snapshots, so their distinction is established by the closedness
diagnostic rather than visual appearance. The confidence bands are
pointwise; values along a trajectory are temporally correlated.

\paragraph{Scope.}
The test uses the absolute annulus complex, with noise ranks $3$, $2$,
and $1$ in runs~1--3. Total injection is deliberately unequal. No separate
coexact-only, relative-boundary, or three-dimensional stochastic run was
performed. Periods were not measured; class preservation is assessed
through closedness and the harmonic projection. The comparisons are with
exact predictions for the fixed implicit scheme, not a time-refinement
or continuum-convergence study.

\section{Discussion and outlook}\label{sec:discussion}
The main distinction is between preserving harmonic coordinates and preserving
a cohomology class. For a fixed $h_0\in\Harm_B$, projected forcing
preserves the affine space $h_0+Q_BH$, whereas exact forcing preserves the
class $h_0+E_B^1$ of closed forms within it. The boundary realization
and Hilbert structure determine these spaces. Deterministic Hodge heat flow
converges to the harmonic projection; continued stochastic forcing can instead
sustain non-harmonic fluctuations.

The deterministic experiments in Section~\ref{sec:num-annulus} reproduce the
predicted boundary-dependent kernel dimensions on an annulus and a solid torus,
with harmonic drift near round-off and non-harmonic decay consistent with the
implicit-Euler spectral factors. The stochastic annulus experiment in
Section~\ref{sec:num-perturbations} compares unprotected, harmonic-projected,
and exact forcing. Among these three constructions, only exact forcing preserves
both harmonic coordinates and closedness. The measured second moments agree
with the discrete predictions; this does not establish continuous-time
covariance accuracy. Time-step and mesh refinement, controlled
projector-error tests, and fixed-weight comparisons remain to be carried out.

Further analysis should address the stability of harmonic projectors and
spectral gaps under perturbations of geometry, weights, and boundary conditions,
with compatible maps identifying the form spaces being compared. Topology
alone does not control the positive spectrum. A filtration introduces a scale
distinct from diffusion time and requires compatible coefficient systems and
maps across scales \cite{WeiWei2025PSL,PersistentHodge2024}. The present
experiments use fixed domains rather than filtrations.

Time-dependent geometry introduces a separate issue. After identifying the
state spaces with a fixed Hilbert space, a deterministic differentiable family
of harmonic projectors $P(t)$ satisfies, for sufficiently regular solutions,
\begin{equation}\label{eq:moving-projector}
\dt(P(t)\omega_t)
=\dot P(t)\omega_t\,\dt t+P(t)\,\dt\omega_t.
\end{equation}
Even when $P(t)G_t=0$ and $P(t)\Delta(t)=0$, the projector-variation term
can change the projected state. Moving-domain models therefore require a
specified transport and control of this additional term.

Beyond constant coefficients, a tractable continuous setting is a finite-rank
local system represented by a flat vector bundle $E\to M$
\cite{HansenGhrist2019}. Its twisted differential satisfies
$d_E^{q+1}d_E^q=0$. With a bundle metric and compatible boundary domains,
one can define
\begin{equation}\label{eq:twisted-hodge-laplacian}
\Delta_{q,E,B}=d_{E,B}^{q-1}(d_{E,B}^{q-1})^*
  +(d_{E,B}^{q})^*d_{E,B}^{q}.
\end{equation}
Under appropriate elliptic and closed-range hypotheses, the kernel represents
absolute or relative cohomology with local coefficients. Extending the
stochastic preservation arguments also requires the corresponding initial-data
and noise-integrability assumptions. More general cellular sheaves with
noninvertible restrictions or varying stalk dimensions do not by themselves
specify a continuum differential complex and boundary realization.

The conditional reverse formulation in Section~\ref{sec:reverse} has not yet
been tested numerically. A finite-dimensional Gaussian-mixture benchmark with
known conditional scores could separate score approximation from terminal-law
and numerical integration errors. Biomedical applications require a specified
image-to-form representation and boundary model, followed by independent
evaluation of predictive performance and decoded-image topology. Preservation
of form coordinates alone does not supply these application-level guarantees.

\section*{Acknowledgment}

The work of CCL was supported in part by Middle Tennessee State University. 
 The work of SC and GWW was supported in part by the University of Georgia, Georgia Research Alliance, and   NIH grant R35GM148196. 

\section*{AI Usage Disclosure}
The authors used Claude (Anthropic) for grammar and \LaTeX{} checking, proofreading,
literature searching, and generating Python code to check the examples.

\end{document}